\documentclass[12pt]{amsart}
\usepackage{amscd}
\def\NZQ{\Bbb}               % the font for N,Z,Q,R,C

\def\ZZ{{\NZQ Z}}

\def\frk{\frak}               % font for "Fraktur"

\def\mm{{\frk m}}

\def\opn#1#2{\def#1{\operatorname{#2}}} % to make
\opn\chara{char}
\opn\length{\ell}
\opn\pd{pd}
\opn\rk{rk}
\opn\projdim{proj\,dim}
\opn\hdeg{hdeg}
\opn\ideg{ideg}
\opn\rank{rank}
\opn\dist{dist}
\opn\depth{depth}
\opn\grade{grade}
\opn\height{height}
\opn\embdim{emb\,dim}
\opn\codim{codim}

\opn\Tr{Tr}
\opn\bigrank{big\,rank}
\opn\superheight{superheight}\opn\lcm{lcm}
\opn\trdeg{tr\,deg}%
\opn\reg{reg}
\opn\lreg{lreg}
\opn\ord{ord}
\opn\div{div}
\opn\Div{Div}
\opn\cl{cl}
\opn\Cl{Cl}
\opn\Spec{Spec}
\opn\Supp{Supp}
\opn\supp{supp}
\opn\Sing{Sing}
\opn\Ass{Ass}
\opn\Ann{Ann}
\opn\Rad{Rad}
\opn\Soc{Soc}
\opn\Ker{Ker}
\opn\Coker{Coker}
\opn\Im{Im}
\opn\Hom{Hom}
\opn\Tor{Tor}
\opn\Ext{Ext}
\opn\End{End}
\opn\Aut{Aut}
\opn\id{id}

\opn\nat{nat}
\opn\pff{pf}%   \pf exists already
\opn\Pf{Pf}
\opn\GL{GL}
\opn\SL{SL}
\opn\mod{mod}
\opn\ord{ord}
\opn\lcm{lcm}
\opn\gcd{gcd}
\opn\Tor{Tor}
\opn\aff{aff}
\opn\con{conv}
\opn\relint{relint}
\opn\st{st}
\opn\lk{lk}
\opn\cn{cn}
\opn\core{core}
\opn\vol{vol}
\opn\link{link}
\opn\star{star}
\opn\gr{gr}

\def\pot#1#2{#1[\kern-0.28ex[#2]\kern-0.28ex]}

\opn\dirlim{\underrightarrow{\lim}}
\opn\inivlim{\underleftarrow{\lim}}
\let\sect=\cap

\let\Dirsum=\bigoplus

\let\mcone= * % Tentative (by TAKAYAMA)

\let\to=\rightarrow

\def\Implies{\ifmmode\Longrightarrow \else
     \unskip${}\Longrightarrow{}$\ignorespaces\fi}
\def\implies{\ifmmode\Rightarrow \else
     \unskip${}\Rightarrow{}$\ignorespaces\fi}
\def\iff{\ifmmode\Longleftrightarrow \else
     \unskip${}\Longleftrightarrow{}$\ignorespaces\fi}

\let\:=\colon
\newtheorem{Theorem}{Theorem}[section]
\newtheorem{Lemma}[Theorem]{Lemma}
\newtheorem{Corollary}[Theorem]{Corollary}
\newtheorem{Proposition}[Theorem]{Proposition}
\newtheorem{Remark}[Theorem]{Remark}

\newtheorem{Definition}[Theorem]{Definition}

\let\epsilon\varepsilon
\let\phi=\varphi
\let\kappa=\varkappa
\opn\inii{in}
\opn\inim{inm}
\opn\set{set}
\def\pnt{{\raise0.5mm\hbox{\large\bf.}}}
\def\lpnt{{\hbox{\large\bf.}}}

\begin{document}

\title{Resolutions of facet ideals}
\author{Xinxian Zheng\\ \\
   Fachbereich Mathematik und Informatik \\
  Universit\"at-GHS Essen, 45117 Essen, Germany\\ 
e-mail: xinxian.zheng@uni-essen.de}

%\date{}

\begin{abstract}
In this paper we study the resolution of a facet ideal associated with a special class of  simplicial complexes
introduced by Faridi.  These simplicial complexes are called trees, 
and  are a generalization (to higher dimensions)  of the concept of a tree in graph theory. We show that the Koszul %%@
homology of the facet ideal $I$ of a tree is generated  by the homology classes of monomial cycles, determine the %%@
projective dimension and the regularity
of $I$ if the tree is 1-dimensional, show that the graded Betti numbers  of $I$ satisfy an alternating sum property if %%@
the tree is connected in codimension 1, and classify all trees whose facet ideal has a linear resolution. 
\end{abstract}

\maketitle

\section*{Introduction}
With a simplicial complex $\Delta$ one can associate two squarefree monomial ideals: the Stanley-Reisner ideal %%@
$I_{\Delta}$ whose generators correspond to the non-faces of $\Delta$,   or the facet ideal $I(\Delta)$ whose %%@
generators correspond to the facets of $\Delta$. The work of Stanley \cite{S}  has demonstrated that there are  deep %%@
relations  between the combinatorial properties of $\Delta$ and the algebraic properties of $I_{\Delta}$.

Facet ideals for graphs have first been considered by Villareal \cite{V1}. In this special case the facet ideal is %%@
called edge ideal, because its generators correspond to the edges of the graph. In his papers \cite{V2} 
and \cite{V3}, Villareal has shown that the edge ideal is the appropriate algebraic object attached to  
a graph.  Among the graphs the trees are the simplest ones. Faridi generalized in \cite {F1} and \cite {F2} 
the definition of tree to simplicial complexes of any dimension, and also introduced facet ideals to study 
trees. 

In the first section of this paper we introduce the basic notions concerning trees, and give a  
characterization of pure trees which are connected in codimension 1. These type of trees play an 
important role in the following sections. 

Our goal here is to study the Koszul cycles of the facet ideal $I\subset R=K[x_1,\ldots, x_n]$ of a tree.  
By this we mean the cycles of the Koszul complex $K_\lpnt(x,R/I)$ of $R/I$ with respect to 
$x_1,\ldots, x_n$. In  Proposition \ref{monomial cycles} we show that the Koszul homology of the facet 
ideal of a tree has a $K$-basis with homology classes of monomial cycles as its elements. In the particular case of a 
$1$-dimensional tree we even show that the Koszul homology of the edge ideal  is generated as a 
$K$-algebra by the homology classes of linear cycles, see Proposition \ref{algebra} . Using this fact, 
in Corollary \ref{reg and pd}, we determine  the regularity and the projective dimension of the facet 
ideal of a $1$-dimensional tree. Furthermore in Theorem \ref{reg} we show that for the facet ideal $I$ of a 
$1$-dimensional tree,  the regularity of $R/I$ is the maximal number $j$, for which 
there exist $j$ edges which are pairwise disconnected.

In the third section, we consider the facet ideal $I$ of a pure tree  and describe the 
linear part of the resolution of $R/I$, see Proposition \ref{linear part}.  We call a tree whose facet 
ideal has a linear resolution a linear tree. In Proposition \ref{linear resolution} we show that a 
tree is a linear tree if and only if the facet ideal of this tree is a linear quotient ideal and we 
classify (Theorem \ref{all linear trees}) all linear trees of a given dimension.  Moreover in 
Corollary \ref{Betti numbers of lqt}, we determine the Betti numbers of the facet ideal of a linear tree.

In the last section, we show that all trees which are connected in codimension $1$ have the alternating 
sum property, meaning that in each linear strand of the resolution of the facet ideal except for the 
lowest one, the alternating sum of the graded Betti numbers is zero, and for the lowest one it is $-1$. 

I would like to thank  Professor J\"urgen Herzog for many helpful comments and  discussions.  

\section{Facet ideals}

In this section we fix the terminology, review some basic properties of graphs and introduce a notion of 
{\em tree} on simplicial complex given by Faridi. As a main result of this section we give a 
characterization of pure trees which are connected in codimension 1. These type of trees play an 
important role in the following sections. 

\begin{Definition}
\label{facet}
{\em A {\em simplicial complex} $\Delta$ over a set of vertices $V=\{v_1,\ldots, v_n\}$ is a collection of subsets of
$V$ with the property that ${v_i}\in\Delta$ for all $i$, and if $F\in\Delta$ then all the subsets of $F$ are
also in $\Delta$ (including the empty set). An element of $\Delta$ is called a $face$ of $\Delta$, and the
{\em dimension} of a face $F$ of $\Delta$ is defined as $|F|-1$, where $|F|$ is the number of vertices of $F$.
In particular, $\dim \emptyset=-1$.  
The faces of dimension 0 and 1 are called {\em vertices} and {\em edges}. The maximal faces of $\Delta$ 
under inclusion are called {\em facets}}.
\end{Definition}

  The dimension of the simplicial complex $\Delta$ is the maximal dimension of its facets, that is to say
\[
\dim\Delta=\max\{\dim F\: F\in\Delta\}. 
\]  
  
We denote the simplicial complex $\Delta$ with the facets $F_1,\ldots, F_q$ by
\[
\Delta=\langle F_1,\ldots, F_q \rangle, 
\] 
and the facet set of $\Delta$ by ${\mathcal F}(\Delta)$.  A simplicial complex $\Delta$ with only one 
facet is called a {\em simplex}, note that $\emptyset$ is also a simplex. A simplicial complex $\Gamma$ is called 
a {\em subcomplex} of $\Delta$ if ${\mathcal F}(\Gamma)\subset {\mathcal F}(\Delta)$.

\begin{Definition}
\label{facet ideal}
{\em Let $\Delta$ be a simplicial complex with vertices  $v_1,\ldots, v_n$. Let $K$ be a field, 
$x_1,\ldots, x_n$ indeterminates, and $R$ the polynomial ring $K[x_1,\ldots, x_n]$. The ideal $I(\Delta)\subset R$ 
generated by the square-free monomials $x_{i_1}\cdots x_{i_s}$, where $\{v_{i_1},\ldots, v_{i_s}\}$ 
is a facet of $\Delta$, is called the {\em facet ideal of $\Delta$}. 
For a $1$-dimensional tree, the facet ideal is called the {\em edge ideal}.} 
\end{Definition}
  
\begin{Definition}
\label{pure}
{\em Let $\Delta$ be a simplicial complex of dimension $d$. Then $\Delta$ is called
 
(a) {\em pure}, if all of its facets have the same dimension;

(b) {\em connected},  if for any two facets $F$ and $G$ there exists a sequence of facets 
$F=F_0,\ldots, F_n=G$, such that $F_i\cap F_{i+1}\neq \emptyset$  for all 
$i=0,\ldots, n-1$; we call this sequence a {\em chain} between $F$ and $G$, and $n$ is called the 
{\em length} of this chain;

(c)  {\em connected in codimension 1}, if for any two facets $F$ and $G$ with $\dim(F)\geq\dim(G)$, there 
exists a chain ${\mathcal C}:\ F=F_0,\ldots, F_n=G$ between $F$ and $G$ such that $\dim(F_i\cap F_{i+1})
=\dim(F_{i+1})-1$ for all $i=0,\ldots, n-1$.}
\end{Definition}

 The chain ${\mathcal C}$ (in Definition \ref{pure} (c)) is called a {\em proper chain}.  
One can see that in a proper chain $\dim F_{i+1}\leq \dim F_i$ for $i=0,\ldots,n-1$.

\begin{Definition}
\label{irredundant chain}
{\em A (proper) chain ${\mathcal C}$ between $F$ and $G$ is called {\em irredundant} if no  subsequence of this chain 
except ${\mathcal C}$ itself is  a (proper) chain between $F$ and $G$.}
\end{Definition}

\begin{Remark}
{\em  Any (proper) chain, after removing suitable facets in it, becomes an irredundant (proper) chain. 
In fact, let ${\mathcal C}$ be a (proper) chain between $F$ and $G$. 
The set of  (proper) subchains of ${\mathcal C}$ is a partially ordered non-empty set. 
The minimal elements in this set are the irredundant (proper) chains between $F$ and $G$. }
\end{Remark}
 
It is clear that an irredundant proper chain need not to be an irredundant chain. For example, 
$F_0=\{a,b,c\}, F_1=\{a,c,d\}, F_2=\{c,d,e\}$ is an irredundant proper chain between $F_0$ and $F_2$, 
but it is not an irredundant chain.

\begin{Lemma}
\label{irredundant proper chain}
Let ${\mathcal C}\: F=F_0, F_1,\ldots,F_n=G$ be a proper chain between $F$ and $G$. If ${\mathcal C}$ is 
irredundant, then $F_j\neq F_k$ for $j\neq k$, and  $F_i\cap F_{i+1}\not\subseteq F_{l}\cap F_i$ for 
$i=1,\ldots,n-1$, and any $l<i$. 
\end{Lemma}
    
\begin{proof} 
Suppose there exists $k>j$ such that $F_j=F_k$, then $F_0,\ldots,F_j, F_{k+1},\ldots, F_n$ is a 
proper subsequence of ${\mathcal C}$ and it is a proper chain between $F$ and $G$, a contradiction.
 
Thus we may now assume $F_j\neq F_k$ for $j\neq k$. Suppose there exists $i\in [n-1]$, such that 
$F_i\cap F_{i+1}\subseteq F_{l}\cap F_i$ for some $l<i$. Then
\[
F_{l}\cap F_{i+1}\supseteq {(F_{l}\cap F_i)\cap (F_i\cap F_{i+1})} =F_i\cap F_{i+1},
\] 
so  
$\dim (F_{l}\cap F_{i+1})\geq\dim (F_i\cap F_{i+1})=\dim F_{i+1}-1$. On the other hand, since 
$F_{l}\neq F_{i+1}$ both are facets, and $\dim F_{i+1}\leq \dim F_{l}$, it follows that 
$\dim (F_{l}\cap F_{i+1})\leq \dim F_{i+1}-1$. Hence 
$\dim (F_{l}\cap F_{i+1})=\dim F_{i+1}-1=\dim (F_i\cap F_{i+1})$, together with 
$F_{l}\cap F_{i+1}\supseteq F_i\cap F_{i+1}$, we have  $F_{l}\cap F_{i+1}=F_i\cap F_{i+1}$. Then 
$F_0,\ldots, F_{l}, F_{i+1},\ldots, F_n$ is a proper subsequence of ${\mathcal C}$, and it is a proper 
chain between $F$ and $G$, a contradiction. 
\end{proof}

Usually a connected graph is called a tree if it has no cycles. Now we define a very special class of trees which play %%@
an important role in Section 2.

\begin{Definition}
\label{bouquet}
{\em A graph $\Gamma$  with vertex set $\{x,y_1\ldots,y_l\}$, $l\geq 1$, and edges $\{x,y_i\}$  for $i=1,\ldots,l$   %%@
is called a {\em bouquet}. We denote this bouquet by $(x;y_1\ldots,y_l)$. The vertex  $x$ is called the {\em root}, %%@
the vertices  $y_i$ the {\em flowers} and the edges $\{x,y_i\}$ the {\em stems} of this bouquet. }
\end{Definition} 

Let $\Delta$ be a tree. If a subgraph $\Gamma$ of $\Delta$ is a bouquet, then we say $\Gamma$ is a bouquet of %%@
$\Delta$.

In \cite{F1}  Faridi introduced the notion of tree  for higher dimensional simplicial  complexes.

\begin{Definition}
\label{leaf}
{\em Let $\Delta$ be a simplicial complex. A facet $F$ of $\Delta$ is called a {\em leaf} if either $F$ 
is the only facet of $\Delta$, or there exists a facet $G$ in $\Delta$, $F\neq G$, such that 
$F\cap H\subseteq F\cap G$ for any facet $H\in\Delta$, $H\neq F$.

We denote the set of all facets $G\in\Delta$ with this property by ${\mathcal U}_{\Delta}(F)$ and call 
it the {\em universal set} of $F$ in $\Delta$. }
\end{Definition}

  For a facet $F$ of $\Delta$, if $x$ is a vertex of $F$ and $x$ does not belong to any other facets of 
$\Delta$, then we call $x$ a {\em free vertex} of $F$ in $\Delta$. It is clear that if $F$ is a leaf of  $\Delta$, 
then $F$ has at least one free vertex. But the converse is not true, even if $\Delta$ is pure.

For example, $\Delta=\langle \{a,b,c\}, \{c,d,e\}, \{e,f,g\} \rangle$ is a pure simplicial complex, the 
facet $\{ c,d,e \}$ has a free vertex $d$, but it is not a leaf.

It is easy to see that $F\in {\mathcal F}(\Delta)$ is a leaf of $\Delta$, if and only if 
$\langle F\rangle  \sect\Gamma$ is a simplex, where 
$\Gamma=\langle {\mathcal F}(\Delta)\setminus \{F\}\rangle$ is the subcomplex of $\Delta$. 
 
\begin{Lemma}
\label{no leaf}
Let ${\mathcal C}\: F_0,\ldots, F_n$ be an irredundant chain in a simplicial complex. Then 
$F_p\sect F_q=\emptyset$ for any $p\in\{0,\ldots,n\}$ and any $q\neq p-1,p,p+1$. Furthermore,  
$F_i$ is not a leaf of $\Gamma =\langle F_0,\ldots, F_n\rangle$ for $i=1,\ldots, n-1$. 
\end{Lemma}

\begin{proof}
Suppose there exists $p\in\{0,\ldots,n\}$ and $q>p+1$ or $q<p-1$, such that $F_p\sect F_q\neq\emptyset$. 
We may assume that $q>p+1$, then $F_0,\ldots,F_p,F_q,\ldots, F_n$ is a chain between $F_0$ and $F_n$, 
a contradiction. 

Suppose $F_j$ is a leaf of $\Gamma$ for some $j\in \{1,\ldots, n-1\}$. Since 
$F_j\sect F_k=\emptyset$ for any $k\neq j-1,j,j+1$, we have $F_{j-1}\sect F_j\subseteq F_j\sect F_{j+1}$ 
or $F_{j}\sect F_{j+1}\subseteq F_{j-1}\sect F_{j}$. We may assume that $F_{j-1}\sect F_j\subseteq F_j\sect F_{j+1}$, 
then $F_j\sect F_{j+1}=(F_{j-1}\sect F_j)\sect (F_j\sect F_{j+1})\subseteq F_{j-1}\sect F_{j+1}$. 
On the other hand, since ${\mathcal C}$ is a chain,  $F_{j}\sect F_{j+1} \neq \emptyset$,   
hence $F_{j-1}\sect F_{j+1}\neq\emptyset$. It follows that
$F_0,\ldots, F_{j-1}, F_{j+1},\ldots, F_n$ is a chain. This contradicts our assumption that 
${\mathcal C}$ is irredundant.
\end{proof} 

We have seen that  an irredundant proper chain need  not to be an irredundant chain. But as im Lemma \ref{no leaf} 
we also have:

\begin{Lemma}
\label{proper no leaf}
Let ${\mathcal C}\: F_0,\ldots, F_n$ be an irredundant proper chain in a simplicial complex, and let 
$\Gamma=\langle F_0,\ldots, F_n\rangle$. Then $F_i$ is not a leaf of $\Gamma$, 
for $i=1,\ldots, n-1$. 
\end{Lemma}

\begin{proof}
Suppose $F_i$ is a leaf of $\Gamma$ for some $i\in [n-1]$. Then there exists an integer 
$k\neq i$ such that $F_i\sect F_{i+1}\subseteq F_i\sect F_k$. Since  ${\mathcal C}$ is an irredundant 
proper chain, it follows from Lemma \ref{irredundant proper chain} that $k>i$. 

For each $k\geq i+1$, we have $\dim(F_i\sect F_{i+1})=\dim F_{i+1}-1\geq\dim F_k-1\geq\dim (F_i\sect F_k)$. 
It follows that  $F_i\sect F_{i+1}=F_i\sect F_k$. 
So $F_0,\ldots,F_i,F_k,\ldots,F_n$ is a proper chain between $F_0$ and $F_n$, a contradiction. 
\end{proof}

\begin{Definition}[Faridi]
\label{tree}
{\em Let $\Delta$ be a connected simplicial complex.  Then $\Delta$ is called a  {\em tree} if every 
nonempty subcomplex   of $\Delta$ has a leaf. A simplicial complex $\Delta$ with the property that every 
connected component is a tree is called a {\em forest}. }
\end{Definition}

As a main result of this section we want to characterize when a pure tree is connected in codimension 1. 
For this purpose we recall the definitions of star and link of a face, see \cite{BH} [Definition 5.3.4].

Let $\Delta$ be a simplicial complex, and $F\in\Delta$. Then {\em star} of $F$ is the set 
\[
\st_{\Delta} F=\{ G\in \Delta \: {F\cup G}\in \Delta\},
\] 
and the {\em link} of $F$ is the set 
\[
\lk_{\Delta} F=\{G\in \Delta \: {F\cup G}\in \Delta,\quad F\cap G=\emptyset\}.
\]

To simplify notation we occasionally omit the index $\Delta$ in $\st_{\Delta}$ or $\lk_{\Delta}$.
Note that $\lk F\subset\st F$, and both are simplicial complex. Furthermore, $\lk_{\Delta} F$ is a 
subcomplex of $\Delta$. Indeed one has ${\mathcal F}(\st G)=\{F\in {\mathcal F}(\Delta)\: G\subset F\}$, and 
${\mathcal F}(\lk  G)=\{F\setminus G\: F\in {\mathcal F}(\st G)\}$. 

We refer the reader to \cite{BH} to see that these notations are crucial in the analysis of the local
cohomology of a Stanley-Reisner ring.

\begin{Proposition}
\label{link}
Suppose that $\Delta$ is a pure tree of dimension $d$. Then the following are equivalent:
\begin{enumerate}
\item[(i)] for all $G\in\Delta$ with $\dim G\leq d-2$, $\lk G$ is connected;
\item[(ii)] $\Delta$ is connected in codimension 1. 
\end{enumerate}
\end{Proposition}

\begin{proof}
(i)$\Rightarrow$(ii): Suppose $\Delta$ is not connected in codimension 1. Then there exists 
$F,H\in{\mathcal F}(\Delta)$ such that there is no proper chain between $F$ and $H$.
Since $\Delta$ is a tree, it is connected, and hence there exists a chain $F=H_0,H_1,\ldots,H_q=H$ 
between $F$ and $G$. Let 
$a=\min\{ \dim(H_i\cap H_{i+1})\: i=0,\ldots,q-1\}$. Since this chain is not proper we have 
$0\leq a<d-1$. We may assume that there is no other chain $F=K_0,\ldots,K_p=H$ in $\Delta$, such that 
$\min \{\dim(K_i\cap K_{i+1})\: i=0,\ldots,p-1\}>a$, otherwise we take this chain instead of 
$H_0,\ldots,H_q$. Let $\{i_1,\ldots,i_m\}\subseteq\{0,\ldots,q\}$ be the subset such that 
$\dim(H_{i_j}\cap H_{i_j+1})=a$, we know that $\{i_1,\ldots,i_m\}\neq\emptyset$. 
By the choice of our chain there must exist $j\in\{1,\ldots,m\}$ such that 
there is no chain $H_{i_{j}}=E_0,E_1,\ldots,E_s=H_{i_{j}+1}$ in $\Delta$ such that 
$\min\{\dim(E_i\cap E_{i+1})\: i=0,\ldots,s-1\}>a$.  

Let $G=H_{i_j}\cap H_{{i_j}+1}$, then $\dim G=a<d-1$. We claim that  
$\lk G$ is not connected. In fact, if $\lk G$ is connected, then there exists a chain $H_{i_j}=D_0, D_1,\ldots,D_l=
H_{{i_j}+1}$ in $\st G$ such that $(D_i\setminus G)\cap (D_{i+1}\setminus G)\neq \emptyset$, 
for any $i=1,\ldots,l-1$. This implies that $\dim(D_i\cap D_{i+1})>a$ for any $i=1,\ldots,l-1$, 
a contradiction to the choice of $j$.

(ii)$\Rightarrow$(i): Suppose there exists  $G\in\Delta$ with $\dim G\leq d-2$, 
such that $\lk G$ is not connected. Then there exist facets $F$ and $H$ in $\st G$ such that there is no 
chain between $F\setminus G$ and $H\setminus G$ in $\lk G$.

Since $\Delta$ is connected in codimension $1$, there exists an irredundant proper chain 
$F=H_0,H_1,\ldots,H_r=H$ between $F$ and $H$. Since $\dim G\leq d-2$, it follows that 
$(H_i\cap H_{i+1})\setminus G\neq\emptyset$, $i=0,\ldots,n-1$. Therefore not all $H_i$ 
belong to $\st G$, because otherwise 
$F\setminus G=H_0\setminus G,H_1\setminus G,\ldots,H_r\setminus G=H\setminus G$ would be a 
chain between  $F\setminus G$ and $H\setminus G$ in $\lk G$.

Let $l=\min \{j\in \{0,\ldots,r\}\: H_{j+1}\notin\st G\}$, and let
$m=\min\{j\in\{l+2,\ldots,n\}\: H_j\\
\in\st G\}$. Now consider the sequence of facets 
$H_l,\ldots,H_m$, it is an irredundant proper chain between $H_l$ and $H_m$, and 
$H_l, H_m\in\st G$, $H_{l+1},\ldots,H_{m-1}\notin\st G$. 

Take the subcomplex 
$\Gamma=\langle H_l,\ldots,H_m\,\rangle$ of $\Delta$. 
Then this subcomplex has no leaf, and so $\Delta$ is not a tree, a contradiction.  
Indeed, since $H_l,\ldots,H_m$ it is an irredundant proper chain, it follows from Lemma \ref{proper no 
leaf} that $H_i$ is not a leaf for $i=l+1,\ldots, m-1$. Now consider the facet $H_l$, and let 
$H_l\cap H_{l+1}=H$.  Then $H$ is a face of $\Delta$ with dimension $d-1$. Let
$\{w\}=H_l\setminus H_{l+1}$ and  $\{u\}=H_{l+1}\setminus H_l$. 
Since $H_{l+1}\notin\st G$, $G\not\subset H$. On the other hand, $H_l\in \st G$, we must have $w\in G$. 
From $H_m\in\st G$ we know $w\in H_m$. That is to say $H_l$ has no free vertex in $\Gamma$, hence $H_1$ 
is not a leaf of $\Gamma$. With the same argument we can show that $H_m$ is not a leaf of $\Gamma$.

\end{proof}

\begin{Corollary}
\label{codimension}
Let $\Delta$ be a pure tree of dimension $d$ and connected in codimension $1$. Then 
 for any facet $F$ of $\Delta$, all the facets of 
      $\langle F\rangle\sect \langle {\mathcal F}(\Delta)\setminus \{F\}\rangle$ are of dimension $d-1$.
\end{Corollary} 

\begin{proof}
 Suppose there exists a facet $F$ of $\Delta$, such that
 $\langle F\rangle\sect \langle {\mathcal F}(\Delta)\setminus \{F\}\rangle$ is not pure of dimension 
$d-1$. Then there exists $H\in {\mathcal F}(\Delta)$ such that $F\cap H=G$ with $\dim G\leq d-2$ and 
$G\not\subset F\cap H'$ for all $H'\in {\mathcal F}(\Delta)\setminus \{F\}$.  

We claim $\lk G$ is not connected. In fact, assume $\lk G$ is connected, then, 
since $F\in {\mathcal F}(\st G)$ if and only if $F\setminus G\in  {\mathcal F}(\lk G)$,   there exists 
a sequence of facets
$F=F_0, F_1,\ldots,F_r=H$ in $\st G$ such that $(F_i\setminus G)\cap (F_{i+1}\setminus G)\neq 
\emptyset$ for
$i=0,\ldots,r-1$. We may assume $F_1\neq F$.  Since $(F\setminus G)\cap (F_1\setminus G)\neq\emptyset$, $G$ is a %%@
proper subset of  $F\cap F_1$,
a contradiction. 

Now Proposition \ref{link} implies $\Delta$ is not connected in
codimension 1, a contradiction to our hypothesis. 
\end{proof}

\begin{Remark}
\label{neq}
{\em Let $\Delta$ be a pure tree of dimension $d$. Even if for any facet $F$ of $\Delta$, all the facets of 
 $\langle F\rangle\sect \langle {\mathcal F}(\Delta)\setminus \{F\}\rangle$ are of dimension $d-1$, 
$\Delta$ may not be connected in codimension $1$.

For example, $\Delta=\langle\{a,b,c\},\{b,c,d\},\{c,e,f\},\{c,f,g\}\rangle$ is pure of dimension $2$, and
for any facet $F$ of $\Delta$, the facets of 
$\langle F\rangle\sect \langle {\mathcal F}(\Delta)\setminus \{F\}\rangle$ are of dimension $1$, but 
$\Delta$ is not connected in codimension $1$.}
\end{Remark}

However we have:

\begin{Corollary}
\label{leaf of original}
Let $\Delta$ be a pure tree of dimension $d$ and connected in codimension $1$, $F$ a facet of $\Delta$.
Then $\Gamma=\langle{\mathcal F}(\Delta)\setminus \{F\}\rangle$ is connected in codimension $1$ if and 
only if $F$ is a leaf of $\Delta$.
\end{Corollary}

\begin{proof}
Assume $F$ is a leaf of $\Delta$. By Lemma \ref{proper no leaf}, for any irredundant proper chain ${\mathcal
C}: F_0,\ldots,F_l$ in $\Delta$, $F\neq F_i$ for any $i\in [l-1]$. Hence $\Gamma$ is connected in codimension
$1$.

Now assume $\Gamma$ is connected in codimension $1$. By Corollary \ref{codimension}, 
$\langle F\rangle\cap\Gamma$ is a pure simplicial complex of dimension $d-1$. Assume $F$ is not a leaf 
of $\Delta$. Then there exist two facets $H_1$ and $H_2$ in $\Gamma$ such that $\dim (F\cap H_i)=d-1$, 
$i=1,2$. Let $G=H_1\cap H_2$. Then $\dim G=d-2$. We may assume $H_1=G\cup\{x_1,x_2\}$ and 
$H_2=G\cup\{x_3,x_4\}$  and $F=G\cup \{x_1,x_4\}$, where $x_i$ are vertices. Since $\Gamma$ is a pure 
tree and connected in codimension $1$, by Proposition \ref{link}, $\lk_{\Gamma}G$ is connected. 
Let $\{x_1,x_2\}=H_1\setminus G=F_1\setminus G,\ldots,F_l\setminus G=H_2\setminus G=\{x_3,x_4\}$ be
an irredundant chain between $\{x_1,x_2\}$ and $\{x_3,x_4\}$ in $\lk_{\Gamma}G$. Then the subcomplex
$\langle F,F_1,\ldots,F_l\rangle$ of $\Delta$ has no leaf, a contradiction. Indeed, each vertex in 
$\langle F,F_1,\ldots,F_l\rangle$ belongs to at least two facets of this subcomplex. 
\end{proof} 

Another consequence of  Corollary \ref{codimension} is

\begin{Proposition}
\label{two leaves}
Let $\Delta$ be a pure tree which is connected in codimension 1, and has more than one facet. Then 
$\Delta$ has at least two leaves.
\end{Proposition}

\begin{proof}
Let $\dim\Delta=d$. Suppose $\Delta$ has only one leaf. Let $F_1$ be this leaf. Since $\Delta$ is connected and 
has more than one facet, there exists a facet $G$ such that $F_1\cap G\neq\emptyset$. Since $\Delta$ 
is pure it follows from  Corollary \ref{codimension} that  there exists a facet $F_2$, such that 
$F_1\cap G\subseteq F_1\cap F_2$ and $\dim (F_1\cap F_2)=d-1$. Let $F_2\setminus F_1=\{ x\}$. Since $F_2$ 
is not a leaf, there exists a facet $H$, such that $\{ x\}\subseteq F_2\cap H$. Again by Corollary 
\ref{codimension} there exists a facet $F_3$, such that $F_2\cap H\subseteq F_2\cap F_3$ and 
$\dim(F_2\cap F_3)=d-1$. It is clear that $F_3\neq F_1$. Since $F_3$ is not a leaf, by the same reason 
there exists a facet $F_4\neq F_2$, and $\dim (F_3\cap F_4)=d-1$, and so on. Since there are only 
finitely many facets, there exist integers $i$ and $j$ with $j<i-1$ such that $F_i=F_j$. If 
$F_{i-1}\cap F_i\neq F_j\cap F_{j+1}$, then the subcomplex $\langle F_j,\ldots,F_{i-1}\rangle$ has no leaf. 
If $F_{i-1}\cap F_i=F_j\cap F_{j+1}$, then the subcomplex $\langle F_{j+1},\ldots,F_{i-1}\rangle$ has no 
leaf. This contradicts our assumption that $\Delta$ is a tree.
\end{proof}

By definition, in a simplicial complex $\Delta$ which is connected in codimension $1$, for any two
facets $F$ and $G$, there exists an irredundant proper chain between $F$ and $G$. For a pure tree we 
even have

\begin {Proposition}
\label {unique chain}
Let $\Delta$ be a pure tree and connected in codimension $1$.  Then for any two facets $F$ and $G$, there exists 
a unique irredundant proper chain between $F$ to $G$.  
\end {Proposition}

\begin {proof}
Suppose ${\mathcal C}\: F=F_0,\ldots ,F_n=G$ and ${\mathcal C '}\: F=G_0,\ldots ,G_m=G$ are two different 
irredundant proper chains between $F$ and $G$. Let $l=\min\{j\: \text{ such that }F_j\neq G_j\}$, and 
$k=\min\{i\: i\geq j+1 \text{ and $F_i=G_t$ for some $t$}\}$. Then, since ${\mathcal C}$ and 
${\mathcal C'}$ both are irredundant, $t>l$ and $G_t\neq G_i$ for any $i\neq t$. Let 
$\Gamma$ be a subcomplex of $\Delta$, such that $F_l,\ldots, F_{k-1}, G_l,\ldots, G_{t-1}\in\Gamma$; if 
$F_{l-1}\cap F_l\neq G_{l-1}\cap G_l$, then let $F_{l-1}\in \Gamma$; if 
$F_{k-1}\cap F_k\neq G_{t-1}\cap G_t$, then let $F_k\in \Gamma$; and there 
are no other facet in $\Gamma$. By Lemma \ref{proper no leaf} one can easily check that $\Gamma$ has no leaf, 
a contradiction since  $\Delta$ is a tree.  
\end {proof}  

According to this proposition, we give the following definition:

\begin{Definition}
\label{distance} 
{\em Let $\Delta$ be a pure tree and connected in codimension $1$. For any two facets $F$ and $G$, the length 
of the unique irredundant proper chain between $F$ and $G$ is called the {\em distance} between $F$ and 
$G$, and denoted by $\dist(F,G)$. 

We call $\max\{\dist(F,G)\: \text{$F$ and $G$ are two facets of }\Delta\}$ the {\em diameter}
of $\Delta$.

If $\Delta$ is a pure forest and each connected component is connected in codimensiom $1$, then for any 
two facets $F$ and $G$ which lie in two different components,  we set $\dist(F,G)=\infty$.}
\end{Definition}

\begin{Remark}
\label{diameter}
{\em Let $\Delta$ be a pure tree and connected in codimension $1$ with diameter $l$, and $F_0,\ldots,F_l$ an
irredundant proper chain of length $l$ in $\Delta$. Then $F_0$ and $F_l$ are leaves of $\Delta$.

Indeed, since $F_l\ldots, F_0$ is also an irredundant proper chain of length $l$, we only need to show $F_0$ is 
a leaf of $\Delta$. Let $d=\dim \Delta$. By Corollary \ref{codimension}, 
$\langle F_0\rangle\sect \langle {\mathcal F}(\Delta)\setminus \{F_0\}\rangle$ is a pure simplicial complex of 
dimension $d-1$. Suppose $F_0$ is not a leaf. Then there exists a facet $F$ of $\Delta$ such that 
$\dim(F_0\cap F)=d-1$ and $F_0\cap F\neq F_0\cap F_1$. Hence $F,F_0\ldots,F_l$ is an irredundant chain 
in $\Delta$ with length $l+1$, a contradiction.}
\end{Remark}

Sometimes we consider a kind of simplicial complex which need not to be a tree, but has some nice properties
like a tree, we call it a quasi-tree.

A connected simplicial complex $\Delta$ is called a {\em quasi-tree}, if there exists an order $F_1,\ldots, F_n$ 
of the facets, such that $F_i$ is a leaf of $\langle F_1,\ldots, F_i \rangle$ for each 
$i=1,\ldots, n$.  
Such an order is called a {\em leaf order}. A simplicial complex $\Delta$ with the property that every 
connected component is a quasi-tree is called a {\em quasi-forest}.

A tree is a quasi-tree, hence for any tree there exists a leaf order of facets. But a quasi-tree need not 
to be a tree. 

For example, $\Delta =\langle \{a,b,c\}, \{b,c,d\}, \{c,d,e\}, \{b,d,f\}\rangle$ is a quasi-tree, but it 
is not a tree, because the subcomplex $\langle\{a,b,c\},\{c,d,e\},\{b,d,f\}\rangle$ has no leaf.
              
\section{On the Koszul cycles of the facet ideal of a tree}

  In the remaining  sections $R=K[x_1,\ldots, x_n]$ denotes the polynomial ring in $n$ indeterminates 
over  the field $K$. Let $M$ be an $R$-module, we denote the Koszul complex $K_\lpnt (x,M)$ of $M$ 
with respect to the sequence $x_1,\ldots, x_n$ by $K_\lpnt (M)$, and for the modules of Koszul cycles, 
Koszul boundaries and the Koszul homology we write  $Z_\lpnt (M)$, $B_\lpnt (M)$ and $H_\lpnt (M)$, 
respectively.

For simplicity, in the remaining sections all simplicial complexes will have the variables  
$x_1,\ldots, x_n$ as vertices.  For a facet $F=\{x_{i_1},\ldots ,x_{i_s}\}$ in $\Delta$, we denote by
$f=x_{i_1}\cdots x_{i_s}$ (in small letter) the monomial in $R$ corresponding to $F$. 

For the proof of the main result of this section which describes the Koszul cycles of certain monomial 
ideals, we need the following general result on the shifts in the resolution of a $\ZZ^n$-graded module:  
Let $M$ be  a finite $\ZZ^n$-graded $R$-module with minimal 
$\ZZ^n$-graded free resolution
\[
\begin{CD}
\cdots @>>> \Dirsum_{a\in \ZZ^n }R(-a)^{b_{1a}} @>>> \Dirsum_{a\in \ZZ^n}R(-a)^{b_{0a}} @>>> M @>>> 0.
\end{CD}
\]
The numbers $b_{ia}$ are called the multigraded Betti numbers of $M$.

We define the {\em support} of an element $a\in \ZZ^n$ to be the set $\supp a=\{i\: a_i\neq 0\}$. Without 
ambiguity, we may set $\supp x^a=\supp a$ for any non-zero monomial. We set $\ZZ^n_+=\{a\in\ZZ^n\: a_i\geq 0 
\text{ for all $i=1,\ldots,n$}\}$. Then we have 

\begin{Lemma}
\label{torsion}
Let $M$ be a torsion-free $\ZZ^n$-graded $R$-module, and $y_1,\ldots, y_s$ a minimal homogeneous 
generating system of $M$. Suppose that $\supp (\deg(y_i))\subseteq\ZZ^n_+$ and 
$t\notin \supp(\deg(y_i))$ for $i=1,\ldots, s$. Then $t\notin \supp(a)$ for all non-zero multigraded 
Betti numbers $b_{ia}$  of $M$.
\end{Lemma}

\begin{proof}
We prove the assertion by induction on $\projdim(M)$. If $\projdim(M)=0$, then the assertion is obvious.  
Now assume $\projdim(M)>0$, and let $F_\lpnt$ be the minimal multigraded free $R$ resolution of $M$, and 
$\epsilon\: F_0\to M$ the augmentation map. 

Obviously $t\notin \supp(a)$  for all $b_{0a}$ which are non-zero. Let $e_1,\ldots, e_s$ be a 
multigraded basis of $F_0$ with $\epsilon(e_i)=y_i$ for $i=1,\ldots,s$,  
and let $z=\sum c_ie_i$ be a homogeneous element in a minimal homogeneous set of generators of $\Ker(\epsilon)$. Then  
$\deg (z)=\deg (c_i)+\deg (e_i)$ for $i=1,\ldots,s$. By assumption we have $[\deg (e_i)]_t=0$ for $i=1,\ldots,s$. 
Suppose $t\in \supp (\deg(z))$, then $[\deg (c_i)]_t>0$ for all $i$ with $c_i\neq 0$. 
This implies that there exist $c_i'\in R$ such that $c_i=x_tc'_i$ for all $i$. So we have $z=x_t\sum c'_ie_i$ and so %%@
$x_t\sum c'_iy_i=0$. Since $M$ is a torsion-free module, it follows that $\sum c'_iy_i=0$, and 
hence $\sum c'_ie_i\in \Ker (\epsilon)$. That is to say,  $z\in \mm \Ker (\epsilon)$, where 
$\mm=(x_1,\ldots,x_n)$, contradicting the assumption  that $z$ belongs to a minimal homogeneous generating system of 
$\Ker (\epsilon)$. Therefore $t$ does not belong to the support of any element in a minimal set of generators of %%@
$\Ker(\epsilon)$. Since $\Ker(\epsilon)$ is torsion free and $\projdim(\Ker(\epsilon))<\projdim(M)$, the lemma follows %%@
from our induction hypothesis. 
\end{proof}
 
Let $J$ be a monomial ideal. As usual we denote by $G(J)$ the unique minimal set of monomial generators of 
$J$. We put $[n]$ to be the set $\{1,\ldots,n\}$.

\begin {Lemma}
\label {support}
  Let $I$ and $J$  be monomial ideals in $R$ with $G(I)=\{f_1,\ldots, f_m\}$ and 
$G(J)=\{f_1,\ldots,f_{m-1}\}$, and let $b$ be the multidegree of $f_m$. If there exists $t\in [n]$, 
such that $t\in \supp(b)$, but $t\notin \supp(\deg(f_i))$ for $i=1,\ldots,m-1$. Then 
$t\in \supp(a)$  for all $a$ with $b_{ia}(R/(J:I)(-b))\neq 0$.
\end {Lemma}

\begin {proof}
  Let $F_\lpnt$ be the minimal $\ZZ^n$-graded free resolution of $R/(J:I)$, then $F_\lpnt (-b)$ is the 
minimal $\ZZ^n$-graded free resolution of $R/(J:I)(-b)$. Since $t\notin \supp(\deg(f_i))$ for 
$1\leq i\leq m-1$, $t$ is not in the support of the  elements of $G(J:I)=G(J:f_m)$. This is because $G(J:f_m)$ 
is a subset of $\{f_1/\gcd(f_1,f_m),\ldots,f_{m-1}/\gcd(f_{m-1},f_m)\}$.
Applying Lemma \ref{torsion} to  $J:I$, we have that $t\notin \supp(a)$  for all 
$b_{ia}(J:I)\neq 0$. Hence $t\notin \supp(a)$ for all $b_{ia}(R/(J:I))\neq 0$. 
Since $b_{ia}(R/(J:I)(-b))=b_{i,a-b}(R/(J:I))$, we have $t\notin\supp(a-b)$, for any $b_{ia}(R/(J:I)(-b))\neq 0$.  
But $t\in \supp(b)$, hence $t\in \supp(a)$.
\end {proof}

\begin{Theorem}
\label{basic}
Let $J\subset R$ be a monomial ideal, $f\in R\setminus J$ a monomial and let $I=(J,f)$. 
Suppose that there exists an integer $t$ such that $x_t$ divides $f$, but $x_t$ does not divide any 
$g\in G(J)$. Then for all $i>0$ there exist short exact sequences
\[
\begin{CD}
0@>>> H_i(R/J)@>>> H_i(R/I)@>\delta >> H_{i-1}(R/(J:I)(-b))@>>> 0, 
\end{CD}
\]
where $b$ is the multidegree of $f$, and for each homology class $[z]\in H_{i-1}(R/(J:I)(-b))$ the 
homology class $[(-1)^{\deg z}z\wedge (f/x_t)e_t]$ is a preimage of $[z]$ under $\delta$. 
\end{Theorem}

\begin{proof}
  From the exact sequence
\begin{eqnarray}
\label{diagram1}
\begin{CD}
0@>>> R/(J:I)(-b)@>>> R/J@>>> R/I@>>> 0,  
\end{CD}
\end{eqnarray}
we get the long exact sequence
\[
\begin{CD}
\cdots @>>> H_i(R/(J:I)(-b))@>>> H_i(R/J)@>>> H_i(R/I)\\
@>\delta >> H_{i-1}(R/(J:I)(-b))@>>> \cdots,  
\end{CD}
\]
  Let $F_\lpnt$ be the minimal $\ZZ^n$-graded free resolution of $R/J$, then 
\[
\Tor^R_i(K,R/J)=\Dirsum K(-a)^{b_{ia}(R/J)}=H_i(R/J),
\] 
and 
\[
\Tor^R_i(K,R/(J:I)(-b))=\Dirsum K(-a)^{b_{ia}(R/(J:I)(-b))}=H_i(R/(J:I)(-b)).
\]
From Lemma \ref{torsion} and Lemma \ref{support} we know that $t\notin \supp(a)$ for all 
$b_{ia}(R/J)\neq 0$, but $t\in \supp(a)$ for all $b_{ia}(R/(J:I)(-b))\neq 0$. Since 
$H_i(R/(J:I)(-b))\to H_i(R/J)$ is a homogeneous homomorphism, it must be the zero map. Hence we have the 
exact sequence as required.
 
To show $[(-1)^{\deg z}z\wedge (f/x_t)e_t]$ is the preimage of $[z]$, we only need to show
\[
d((-1)^{\deg z}z\wedge (f/x_t)e_t)=fz\quad \text{in}\quad K_{\lpnt}(R/J).
\]

In fact, $d((-1)^{\deg z}z\wedge (f/x_t)e_t)=(-1)^{\deg (z)}d(z)\wedge ((f/x_t)e_t)+fz$. 
Now since  $z\in Z_{i-1}(R/(J:I)(-b))$, it follows that $d(z)\in (J:I)K_{i-2}(R)$, and hence 
$x_t(f/x_t)d(z)=fd(z)\in JK_{i-2}(R)$. Since $x_t$ does not divide any $g\in G(J)$, we have $J=J:x_t$, 
and so $(f/x_t)d(z)\in JK_{i-2}(R)$.  Hence $d(z)\wedge (f/x_t)e_t\in JK_{i-1}(R)$. That is to say, 
$d(z)\wedge (f/x_t)e_t=0$ in $K_{i-1}(R/J)$.
\end{proof}

\begin{Corollary}
\label {regular sequence}
Let $L\subset R$ be a graded ideal, and $x_{i_1},\ldots ,x_{i_s}$ a regular sequence on $R/L$. If $B$ is 
a $K$-basis of $H_\lpnt (R/L)$, then $\{ [z\wedge e_I]\: [z]\in B, \quad I\subset \{i_1,\ldots ,i_s\} \}$ 
is a $K$-basis of $H_\lpnt (R/L+(x_{i_1},\ldots ,x_{i_s}))$.
\end{Corollary}

\begin{proof}
We may assume that $s=1$. The general case is done by induction on $s$. Since $x_{i_1}$ is regular on $R/L$, $x_{i_1}$ %%@
does not divide any $g\in G(L)$. Therefore the result  follows from Theorem \ref{basic}. 
\end{proof}

\begin{Corollary}
\label{injective}
Let $I$, $I'$ be monomial ideals in $R$ with $G(I)=\{f_1,\ldots,f_m,\ldots,f_l\}$ and 
$G(I')=\{f_1,\ldots,f_m\}$. If for any $i\geq m+1$ there exists a variable which divides $f_i$ 
but  does not divide $f_j$ for any $j<i$. Then the map $H_i(R/I')\to H_i(R/I)$ is injective.
\end{Corollary}

\begin{proof}
The statement follows immediately from Theorem \ref{basic} by induction on \\
$l-m$.
\end{proof}

\begin {Corollary}
\label {Betti number}
Let $R$, $I$, $J$ be as in Theorem {\em \ref {basic}}. Then  we have 
\[
b_{ia}(R/I)=b_{ia}(R/J)+b_{i-1,a-b}(R/(J:I))\quad \text{for all $i>0$ and $a\in\ZZ^n$}.
\]
\end {Corollary}

For another main result of this section, we need the following concept:

\begin{Definition}
\label{def monomial cycles} 
{\em Let $I$ be a monomial  ideal of $R$. A cycle $z$ of $K_\lpnt (R/I)$ is called a monomial cycle if there 
exists $L\subset [n]$ and a monomial $f$, such that $z=fe_L$.}  
\end{Definition}

Even if $I$ is a square-free monomial ideal, $H_\lpnt(R/I)$ may not 
be generated by homology classes of monomial cycles. For example, let $R=K[x_1,x_2,x_3,x_4]$ and
$I=(x_1x_2,x_2x_3,x_3x_4,x_4x_1)$. Then $z=x_1e_2\wedge e_3\wedge e_4+x_3e_1\wedge e_2\wedge e_4$ is a 
cycle, but $z$ is not homologous to a monomial cycle. In fact, a boundary $b\in B_3(R/I)$ is of the form 
$d(fe_1\wedge e_2\wedge e_3\wedge e_4)$. So $z$ can not be a monomial cycle.

However for the facet ideal $I$ of a forest, we have $H_\lpnt (R/I)$ is generated by homology classes of 
monomial cycles. To prove this we need the following lemma.

\begin{Lemma}
\label{quotient ideal}
Let $\Delta$ be a forest and  $I$ its facet ideal. If $F$ is any facet of $\Delta$ and $J$ is the 
ideal generated by $G(I)\setminus \{f\}$. Then the simplicial complex $\Delta'$ with facet ideal $J:I$ 
is again a forest. 
\end{Lemma}

\begin{proof}
Note that ${\mathcal F}(\Delta')$ is a subset of $\{G\setminus F\: G\in {\mathcal F}(\Delta)\}$. Suppose $\Delta'$ is %%@
not a forest. Then there exist facets $F_1,\ldots,F_p$ of $\Delta$, such that the
subcomplex $\langle F_1\setminus F,\ldots,F_p\setminus F\rangle$ of $\Delta'$ has no leaf. Since $\Delta$ 
is a forest, the subcomplex $\langle F_1,\ldots,F_p\rangle$ has a leaf $F_i$. Hence there exists a 
integer $k\in\{1,\ldots,p\}$ and $k\neq i$, such that $F_j\sect F_i\subseteq F_k\sect F_i$ for any 
$j\neq i$. Therefore $(F_j\sect F_i)\setminus F\subseteq (F_k\sect F_i)\setminus F$ for any $j\neq i$, 
and hence $(F_j\setminus F)\sect (F_i\setminus F)\subseteq (F_k\setminus F)\sect (F_i\setminus F)$ for any 
$j\neq i$. So $F_i\setminus F$ is a leaf of $\langle F_1\setminus F,\ldots,F_p\setminus F\rangle$, a
contradiction.
\end{proof}

\begin{Proposition}
\label{monomial cycles}
Let $\Delta$ be a forest and $I$ its facet ideal. Then $H_r(R/I)$ has the $K$-basis 
\[
{\mathcal M_r}=\{[fe_{i_1}\wedge\cdots\wedge e_{i_r}]\: fe_{i_1}\wedge\cdots\wedge e_{i_r}\text{ is a monomial
cycle.}\}
\]
\end{Proposition}

\begin{proof}
Let $\Delta =\langle F_1\ldots,F_m\rangle$ where  $F_1\ldots,F_m$ is a leaf order. We prove the assertion 
by induction on $m$. The case $m=1$ is trivial. Since $F_m$ is a leaf we may assume that $f_m=h x_t$, 
where $h$ is a monomial and $x_t\in {F_m\setminus {\bigcup _{j<m}F_j}}$. By Theorem \ref{basic}, we have 
short exact sequences
\[
\begin{CD}
0@>>> H_r(R/J)@>>> H_r(R/I)@>>> H_{r-1}(R/(J:I)(-b))@>>> 0, 
\end{CD}
\]
where $J=(f_1,\ldots, f_{m-1})$ and $b$ is the multidegree of $f_m$. By Lemma \ref{quotient ideal}, $J:I$ 
is a facet ideal of a forest and it has at most $m-1$ facets. Again use Theorem \ref{basic} we have $[z]$,
$[z'\wedge(f_m/x_t)e_t]$ are basis elements of $H_r(R/I)$, where $[z]$ and $[z']$ are basis elements of 
$H_r(R/J)$ and $H_{r-1}(R/(J:I)(-b))$, respectively. And by induction hypothesis $z$ and $z'$ can be 
choosen as monomial cycles. 
\end{proof}

\begin{Definition}
\label{linear cycle}
{\em Let $I$ be a monomial ideal and let $d$ be the least degree of its generators. A monomial cycle $z=fe_L$ 
in $K\lpnt (R/I)$ is called  {\em linear} if $f$ is a monomial of degree $d-1$.}
\end{Definition}
 
\begin{Remark}
\label{k-basis}
{\em Let $\Delta$ be a 1-dimensional forest with edge ideal $I$.  Then the linear monomial cycles 
are of the form 
\[
x_le_{l_1}\wedge \cdots \wedge e_{l_r},
\] 
where $\{x_l, x_{l_i}\}$  is an edge of $\Delta$,  $i=1,\ldots, r$. Hence it follows from Proposition \ref{monomial %%@
cycles}  that  the set 
\[
{\mathcal B}_r =\{[z(b)]\: b=(x_l;x_{l_1},\ldots,x_{l_r})\text{ is a bouquet of }\Delta\text{ with $r$ flowers}\}
\] 
is a $K$-basis of $H_r(R/I)_{r+1}$, where $z(b)=x_le_{l_1}\wedge \cdots \wedge e_{l_r}$.

}
\end{Remark}

\begin{Proposition}
\label{algebra}
Let $\Delta$ be a forest of dimension $1$ and $I$ its facet ideal. Then as a $K$-algebra, $H_\lpnt(R/I)$ 
is generated by the homology classes of linear monomial cycles. 
\end{Proposition}

\begin{proof}
Let $fe_L$ be an arbitrary monomial cycle, and let $i\in L$. Then $fx_i\in I$, and hence there exists a 
generator $f_1\in G(I)$ such that $fx_i=f_1g$. Since $f\not\in I$, we conclude that $x_i$ divides $f_1$. 
Then $f=(f_1/x_i)g$. Now let $L_1=\{l\in L\: (f_1/x_i)x_l\in I\}$, and $L_2=L\setminus L_1$. Note that 
$i\in L_1$ and that $fe_L=(f_1/x_i)e_{L_1}\wedge ge_{L_2}$, where $(f_1/x_i)e_{L_1}$ is a linear cycle.  
If $g=1$, then $fe_L$ is a linear cycle, and if $g\neq 1$ but $L_2=\emptyset$, then $fe_L$ is a boundary. 
Thus we may assume that $g\neq 1$ and $L_2\neq \emptyset$, and have to show that $ge_{L_2}$ is a cycle. 
Then we can proceed by induction on the degree of $f$.

Suppose $gx_s\not \in I$ for some $s\in L_2$. Since $f=(f_1/x_i)g$ we have  $((f_1/x_i)g)x_s\in I$. Let 
$f_1/x_i=x_r$. By the choice of $L_2$ it follows that $x_rx_s\not\in I$. Therefore there must exist $x_t$ 
dividing $g$ such that $x_tx_s\in I$. This implies $gx_s\in I$, a contradiction.  
\end{proof}

%For any forest and quasi-forest $\Delta$, we can give the facets a leaf order, such that it satisfies the
%condition in Proposition \ref{monomial cycles}, so in these cases the Koszul homology $H(R/I)$ for a 
%facet ideal  $I$ has the nice property of Proposition \ref{algebra}.

\begin{Corollary}
\label{reg and pd}
Let $\Delta$ be a 1-dimensional forest with edge ideal $I$. Then
\begin{enumerate}
\item[(i)] $\reg (R/I)$ is the maximal number $j$ for which there exist linear monomial cycles $z_i$ such that 
$[z_1]\cdots [z_j]\neq 0$;
\item[(ii)] $\pd (R/I)$ is the maximum  among the sums $\sum_{i=1}^j k_i$ for which there exist linear 
cycles $z_i\in Z_{k_i}(R/I)$ such that $[z_1]\cdots [z_j]\neq 0$.   
\end{enumerate}
\end{Corollary}

\begin{Proposition}
\label{homology}
Let $\Delta$ be a $1$-dimensional forest with edge ideal $I$, and let
$b_1=(x_{q_1};x_{q_{11}},\ldots,x_{q_{1p_1}}),\ldots,b_l=(x_{q_l};x_{q_{l1}},\ldots,x_{q_{lp_l}})$ be bouquets in %%@
$\Delta$. Then the
following statements are equivalent:
\begin{enumerate}
\item[(i)] $[z(b_1)]\cdots [z(b_l)]\neq 0$.
\item[(ii)] The set of bouquets $b_1,\ldots,b_l$ satisfies the following conditions:
   \begin{enumerate}
   \item[(a)] All vertices occurring in these bouquets are pairwise distinct.
   \item[(b)] The roots of any two bouquets have no common edge.
   \item[(c)] For all bouquets $b_i$ there exists at least one flower which has no 
              common edge with the root of $b_j$ for all $j\neq i$.
   \end{enumerate} 
\end{enumerate}
\end{Proposition}

\begin{proof}
(i)\implies (ii): It is clear that if (a) or (b) not holds, then $[z(b_1)]\cdots[z(b_l)]=0$. Suppose 
there exists an $i$, such that each flower of $b_i$ has common edge with the root of some $b_j$. Since 
$d(e_{q_i}\wedge e_{q_{i1}}\wedge\cdots\wedge e_{q_{ip_i}})=x_{q_i} e_{q_{i1}}\wedge\cdots\wedge 
e_{q_{ip_i}}-x_{q_{i1}} e_{q_i}\wedge e_{q_{i2}}\wedge\cdots\wedge 
e_{q_{ip_i}}+ \cdots +(-1)^{p_i} x_{q_{ip_i}} e_{q_i}\wedge e_{q_{i1}}\wedge\cdots\wedge e_{q_{i,p_i-1}}$, 
we have 
\[
[z(b_i)]=\sum_{k=1}^{p_i} [(-1)^{k+1} x_{q_{ik}} e_{q_i}\wedge e_{q_{i1}}\wedge\cdots\wedge\widehat 
{e_{q_{ik}}}\wedge\cdots\wedge e_{q_{ip_i}}].
\]
Since $x_{q_{ik}}$ has a common edge with the root of some $b_j$ for all $k\in\{1,\ldots,p_i\}$, 
we have  $[z(b_i)][z(b_j)]=0$, a contradiction.

(ii)\implies (i): We prove the assertion by induction on $l$. The case $l=1$ follows from 
Remark \ref{k-basis}. Let $\Delta '$ be the subforest of $\Delta$ obtained as follows:
If one stem  of our bouquets is a leaf of $\Delta$, then let $\Delta'=\Delta$.  Otherwise let $F_1$ be 
any leaf of $\Delta$, and let $\Delta_1=\langle {\mathcal F}(\Delta) \setminus \{F_1\}\rangle$. Notice 
that $\Delta_1$ is again a forest containing all our bouquets. 
If one stem of our bouquets  is a leaf of $\Delta_1$,  then let $\Delta'=\Delta_1$. Otherwise let $F_2$ 
be any leaf of $\Delta_1$, and let $\Delta_2=\langle {\mathcal F}(\Delta_1) \setminus \{F_2\}\rangle$.
Proceeding in this way we obtain a subforest $\Delta'$ of $\Delta$ such that 
 \[
 \Delta=\langle {\mathcal F}(\Delta'), F_s,\ldots,F_1\rangle, 
 \]
 where $F_r$ is a leaf of $\langle {\mathcal F}(\Delta'), F_s,\ldots,F_r\rangle$ for $r\in\{1,\ldots,s\}$, 
and such that some stem of our bouquets, say $\{x_{q_i},x_{q_{ik}}\}$,  is a leaf of $\Delta'$. 
Let $I'$ be the edge ideal of $\Delta '$,  
$\Gamma=\langle {\mathcal F}(\Delta')\setminus\{x_{q_i},x_{q_{ik}}\}\rangle$  
with edge ideal $J'$, and let $\Gamma '$ be the simplicial complex with facet ideal $J':I'$. 
By Lemma \ref{quotient ideal}, $\Gamma '$ is a forest. 

If $p_i>1$, then $x_{q_{ik}}$ must be the free vertex of $\{x_{q_i},x_{q_{ik}}\}$ in $\Delta '$. 
If $p_i=1$, then $b_i=(x_{q_i};x_{q_{i1}})$. It may be that $x_{q_{i1}}$ is not the free vertex of 
$\{x_{q_i},x_{q_{i1}}\}$ in $\Delta'$. Then we replace $b_i$ by the bouquet $b_i'=(x_{q_{i1}};x_{q_i})$. 

Notice that the bouquets $b_1,\ldots, b_{i-1}, b_i',b_{i+1},\ldots, b_l$ again satisfy all conditions in 
(ii), and since $[b_i]=[b_i']$ we also have  $[b_1]\cdots [b_l]=[b_1]\cdots [b_{i-1}][b_i'][b_{i+1}]\cdots [b_l]$. 
Therefore we may as well assume that in any case the flower $x_{q_{ik}}$ is the free vertex of 
$\{x_{q_i},x_{q_{ik}}\}$ in $\Delta'$. 

It follows from the definition of $\Gamma '$ that all the other flowers of $b_i$ are isolated vertices 
of $\Gamma '$. Recall that a vertex  in a simplicial complex  $\Sigma$ is called isolated if it has no  
common edge with any other vertex  in $\Sigma$.

We distinguish two cases:

Case 1: The root of the bouquet $b_i$ has no common edge with any flower in the other bouquets. 

In this case $b_1,\ldots,b_{i-1},b_{i+1},\ldots,b_l$ are bouquets in $\Gamma '$, and this set of bouquets 
satisfies all conditions in (ii). 
By induction hypothesis, $[b_1]\cdots[b_{i-1}][b_{i+1}]\cdots [b_l]\neq 0$ in $H_\lpnt(R/(J':I'))$. 
Since $x_{q_{im}}$ is an isolated vertices in $\Gamma '$ for $m\in\{1,\ldots,p_i\}$ and $m\neq k$, by 
Corollary \ref{regular sequence}, we have $[b_1]\cdots [b_{i-1}][b_{i+1}]\cdots [b_l]
[e_{q_{i1}}\wedge\cdots \wedge \widehat {e_{q_{ik}}}\wedge \cdots\wedge e_{q_{ip_i}}]\neq 0$ in 
$H_\lpnt(R/(J':I'))$. By Theorem \ref{basic}, for any basis element $[z]$ of $H_{r-1}(R/(J':I')(-2))$, 
$[z\wedge x_{q_i}e_{q_{ik}}]$ is a basis element of $H_r(R/I')$. Since  
$z=z(b_1)\cdots \widehat{z(b_i)}\cdots z(b_l)(e_{q_{i1}}\wedge\cdots \wedge\widehat {e_{q_{ik}}}\wedge
\cdots\wedge e_{q_{ip_i}})$ is a cycle in $K_\lpnt(R/(J':I')(-2))$, it follows that 
$0\neq [z\wedge x_{q_i}e_{q_{ik}}]=[b_1]\cdots [b_l]$ in $H_\lpnt(R/I')$. By Corollary \ref{injective}, 
we have $[b_1]\cdots [b_l]\neq 0$ in $H_\lpnt(R/I)$.

Case 2. There exists an integer $j\neq i$ such that the root $x_{q_i}$ has a common edge with some flower 
of $b_j$.

Let $C$ be the set of integers having this property, and let $j\in C$. Since $\Delta$ is a tree, there 
exists only one flower of $b_j$ which has a common edge with $x_{q_i}$, because otherwise $\Delta$ would 
have a cycle. And by the condition (c) in (ii), we have $p_j>1$. For $j\in C$, let 
\[
 b_j'= \left\{\begin{array}{ll}
               b_j,  &\text{if $j\not\in C$,}\\
               (x_{q_j};x_{q_{j1}},\ldots,\widehat {x_{q_{jk}}},\ldots,x_{q_{jp_j}}), 
                                                       &\text{if $j\in C$ and $\{x_{q_i},x_{q_{jk}}\}$ 
                                                               is an edge.}
                                                                                 
                           \end{array}
                           \right.
\] 
Then  $b_1',\ldots,b_{i-1}',b_{i+1}',\ldots,b_l'$ are bouquets of $\Gamma '$, and this set of bouquets 
satisfies all the conditions in (ii).   For all $j\in C$, let $\{x_{q_i},x_{q_{j_k}}\}$  be the unique 
common edge of the root $x_{q_{i}}$ of $b_i$ with the flower $x_{q_{j_k}}$ in $b_j$. Then $x_{q_{j_k}}$ is an 
isolated vertex of $\Gamma '$. Hence in $\Gamma'$ we are in the same situation as in Case 1, and so as 
before the result follows by induction.
\end{proof}

\begin{Definition}
\label{pairwise disconnected}
{\em Let $\Delta$ be a simple graph, that is, for each edge $\{x_i,x_j\}$ of $\Delta$, 
$x_i\neq x_j$. Two edges $\{x_i,x_j\}$ and $\{x_k,x_l\}$ are called {\em disconnected} if
\begin{enumerate}
   \item[(a)] $\{x_i,x_j\}\cap \{x_k,x_l\}=\emptyset$;
   \item[(b)] $\{x_i,x_k\}$, $\{x_i,x_l\}$, $\{x_j,x_k\}$, $\{x_j,x_l\}$ are not edges of $\Delta$.
   \end{enumerate} } 
\end{Definition}

\begin{Corollary}
\label{product}
Let $\Delta$ be a 1-dimensional forest with edge ideal $I$, and 
$\{x_{i_1},x_{j_1}\},\\
\ldots,\{x_{i_m},x_{j_m}\}$ edges of $\Delta$. Then the following are equivalent:
\begin{enumerate}
\item[(i)] $[x_{i_1}e_{j_1}]\cdots[x_{i_m}e_{j_m}]\neq 0$.
\item[(ii)] The edges $\{x_{i_1},x_{j_1}\},\ldots,\{x_{i_m},x_{j_m}\}$ are pairwise disconnected.
\end{enumerate}
\end{Corollary}

\begin{proof}
Let $b_l=\{x_{i_l};x_{j_l}\}$, $l=1,\ldots,m$. Then $b_l$ is a bouquet with one flower. Notice that  
$b_l'=\{x_{j_l};x_{i_l}\}$ is also a bouquet with one flower of $\Delta$. Since 
$[z(b_l)]=[z(b_l')]$, we have $[z(b_1)]\cdots[z(b_m)]\neq 0$ if and only if 
$[z(b_1)]\cdots[z(b_{l-1})][z(b_l')][z(b_{l+1})][z(b_m)]\neq 0$. 
Hence we may choose  $x_{i_l}$ or $x_{j_l}$ as the root of $b_l$.

(i)\implies (ii):    If $[x_{i_1}e_{j_1}]\cdots[x_{i_m}e_{j_m}]\neq 0$, 
then all conditions in (ii) of Proposition \ref{homology} hold.  Hence all vertices occurring in these 
edges are pairwise distinct, and $\{x_{i_l},x_{j_l}\}$, $l=1,\ldots,m$ are the only edges in the subgraph 
of $\Delta$ restricted to the vertices $\{x_{i_1},\ldots, x_{i_m},x_{j_1}\ldots, x_{j_m}\}$. It follows 
that $\{x_{i_1},x_{j_1}\},\ldots,\{x_{i_m},x_{j_m}\}$ are pairwise disconnected.

(ii)\implies (i): If $\{x_{i_1},x_{j_1}\},\ldots,\{x_{i_m},x_{j_m}\}$ are pairwise disconnected, then the 
set of bouquets $b_1,\ldots,b_m$ satisfies all conditions in (ii) of Proposition \ref{homology}. Hence 
$[x_{i_1}e_{j_1}]\cdots[x_{i_m}e_{j_m}]\neq 0$.
\end{proof}

Moreover, we have

\begin{Corollary}
\label{pick up}
Let $\Delta$ be a $1$-dimensional tree, and $b_1,\ldots,b_l$ bouquets of $\Delta$. If the set of these 
bouquets satisfies the condition {\em (ii)} of Proposition {\em \ref{homology}}, then there exists one 
stem in each bouquet, such that these stems are pairwise disconnected.
\end{Corollary}

\begin{proof}
We refer to the notation in the proof of Proposition \ref{homology}.  By the proof \ref{homology} 
(ii) \implies (i) in each step we get a leaf $\{x_{p_i},x_{p_{ik}}\}$ in the subforest of the previous one. 
The arguments in the proof show that these stems are pairwise disconnected.
\end{proof}

By using Proposition \ref{homology}, Corollary \ref{product} and Corollray \ref{pick up}, we conclude:

\begin{Theorem}
\label{reg}
Let $\Delta$ be a $1$-dimensional forest, $I$ its edge ideal. Then the regularity of $R/I$ is the maximal
number $j$, for which there exist $j$ edges which are pairwise disconnected.
\end{Theorem}

\begin{Remark}
\label{without cycle}
{\em In Theorem \ref{reg}, the assumption that $\Delta$ is a forest  is important. If $\Delta$ has a cycle, 
then the assertion might not be true.

For example, let $\Delta$ be a graph with edge ideal $I=(ab,bc,cd,de,ea)$. Then the regularity of $R/I$ is 
$2$, but the maximal number of the pairwise disconnected  edges in $\Delta$ is $1$.}
\end{Remark} 

\section{Linear trees}

In general, it is not easy to determine the Betti numbers of an $R$-module $M$, but for a facet ideal $I$ 
of a pure tree which is connected in codimension 1, we can describe  the linear part of the resolution of 
$R/I$.

  We know that if $M$ is a graded $R$-module, $z\in R$ is a homogeneous element of degree 1, and $z$ is 
a non-zero divisor of $M$, then 
\begin{eqnarray}
\label{**}
b_{ij}(M/zM)=b_{ij}(M)+b_{i-1,j}(M(-1))=b_{ij}(M)+b_{i-1,j-1}(M). 
\end{eqnarray}
In fact, if $F\lpnt$ is a graded minimal free resolution of $M$, then  the mapping cone of 
$F\lpnt(-1)\stackrel{z}{\longrightarrow} F\lpnt$ is the 
minimal graded free resolution of $M/zM$.

\begin {Lemma}
\label {b_ii}
Let $L$ be a monomial ideal in $R$ with  $G(L)=\{g_1,\ldots,g_l\}$. Suppose that $\deg (g_r)=1$ for 
$r=1,\ldots ,s$, and $\deg(g_r)>1$ for $r=s+1,\ldots ,l$. Then $b_{ii}(R/L)={s\choose i}$.
\end {Lemma}

\begin {proof}
We may assume that $g_i=x_i$ for $i\in[s]$. Then for all $i\in[s]$,  $x_i$ does not divide any $g_j$ for 
$j>s$, because $\{g_1,\ldots, g_l\}$ is a minimal set of generators of $L$. Hence $g_1,\ldots, g_s$ is a  
regular sequence modulo $(g_{s+1},\ldots, g_l)$.  Hence the assertion follows by induction on $s$ from 
(\ref{**}).
\end {proof}

\begin{Definition}
\label{adjacent face}
{\em Let $\Delta$ be a $d$-dimensional pure tree connected in codimension 1 and $G$ a face of dimension $d-1$.
If $G$ is contained in at least two facets of $\Delta$, then we call $G$ an {\em adjacent face}.}
\end{Definition}

\begin {Proposition}
\label {linear part}
Let $\Delta$ be a $d$-dimensional pure tree with $m$ facets, $I$  its facet ideal. Suppose $\Delta$ is 
connected in codimension 1. For each adjacent face $G\in \Delta$, 
let $m(G)=|\{F\in{\mathcal F}(\Delta)\: G\subset F\}|$. 
Then 
\[
 b_{i,i+d}(R/I)= \left\{\begin{array}{ll}
                   m,  &\text{if $i=1$,}\\
               \sum_{G}\binom{m(G)}{i}, &\text{if $i\geq 2$.}
                \end{array}
                           \right.
\] 
\end {Proposition}

\begin {proof}
Let $\Delta =\langle F_1,\ldots ,F_m\rangle$ such that $F_1,\ldots, F_m$ is a leaf order.  We prove the proposition by 
induction on $m$. The case $m=1$ is trivial.  Let $\Gamma=\langle F_1,\ldots ,F_{m-1}\rangle$ and $J$ be the facet %%@
ideal of $\Gamma$. By 
Corollary  \ref{Betti number} and Lemma \ref{b_ii}, we know 
\begin{eqnarray*}
b_{i,i+d}(R/I)&=&b_{i,i+d}(R/J)+b_{i-1,i+d}(R/(J:I)(-(d+1)))\\
      &=&b_{i,i+d}(R/J)+b_{i-1,i-1}(R/(J:I)),
\end{eqnarray*}
and 
\[
b_{i-1,i-1}(R/(J:I))={s\choose {i-1}},
\] 
where $s=|\{F_j\: \dim(F_j\cap F_m)=d-1,\quad j=1,\ldots,m-1\}|$, because $x_i\in J:I$ if and only if $x_i\in
F_j\setminus F_m$ for some $F_j$ in this set.

Let $G$ be an adjacent face of $\Gamma$ and $m'(G)=|\{F\in{\mathcal F}(\Gamma)\: G\subset F\}|$. By our 
induction hypothesis
\[
 b_{i,i+d}(R/J)= \left\{\begin{array}{ll}
               m-1,  &\text{if $i=1$,}\\
               \sum_{G}\binom{m'(G)}{i}, &\text{if $i\geq 2$.}
                           \end{array}
                           \right.
\]

So 
\[
 b_{i,i+d}(R/I)= \left\{\begin{array}{ll}
               m-1+1=m,  &\text{if $i=1$,}\\
               \sum_{G}\binom{m'(G)}{i}+\binom{s}{i-1}=\sum_{G}\binom{m(G)}{i}, &\text{if $i\geq 2$.}
                           \end{array}
                           \right.
\]
\end {proof}

For a $d$-dimensional pure tree $\Delta$, we assign to each face $G$ with dimension $d-1$ an {\em degree},  
namely
\[
\deg (G)=|\{F\: F \text{ is a facet of } \Delta, \text{ such that }G\subset F\}|.
\]
By Proposition \ref{linear part}, $b_{i,i+d}(R/I)=\sum_{G}{\deg (G)\choose i}$ for $i\geq 2$, where $I$ 
is the facet ideal of $\Delta$. (Notice that if $G$ is not an adjacent face, then ${\deg (G)\choose i}=0$ 
for $i\geq 2$.) If $d=1$, then the face of dimension $d-1$ is just a vertex. 

\begin{Remark}
\label{the second betti number}
{\em Let $\Delta$ be a $1$-dimensional tree. Then $b_{2,3}=\sum_v{\deg (v)\choose 2}$, where $v$ runs through all
the vertices of $\Delta$. In \cite {EV}, Eliahou and Villarreal proved that for any graph $G$,
$b_{2,3}=|E(L(G))|-N_t$, where $N_t$ is the number of triangles of $G$ and $L(G)$ is the line graph of $G$. In
the case $G$ is a tree, $N_t=0$, and $|E(L(G))|$ is just $\sum_v{\deg (v)\choose 2}$.}
\end{Remark}

\begin {Lemma}
\label {numbers}
Let $\Delta$ be a $d$-dimensional pure tree and connected in codimension $1$, $V$  the set of faces of
dimension $d-1$, and $O=\sum_{G \in V}\deg(G)$.  Then we have $|{\mathcal F}(\Delta)|-1=O-|V|$.
\end {Lemma}

\begin {proof}
The lemma follows by induction on the number of facets, observing  that when we add a leaf to the tree, $O$ 
will increase by $d+1$, and $|V|$  by $d$. 
\end {proof}

For a $d$-dimensional pure tree $\Delta$ which is connected in codimension $1$, let $b_0'=|V|$, $b_1'=O$, 
and $b_i'=b_{i,i+d}$ for $i\geq 2$. 
By using the well-known binomial formula $\sum_{i=0}^{n} {(-1)^i{n\choose i}}=0$, one  sees
that $\sum_i (-1)^ib_i'=0$. Hence together with Lemma \ref {numbers} we have

\begin {Proposition}
\label {alternating sum}
Let $\Delta$ be a $d$-dimensional pure tree with the facet ideal $I$. Suppose $\Delta$ is connected in 
codimension $1$.  Then
\[ 
1+\sum_{i>0} (-1)^ib_{i,i+d}=0.
\]
\end {Proposition} 

In the next section, we will have another property on the Betti numbers of facet ideals which generalizes
this proposition.

\begin{Definition}
\label{linear q-ideal}
{\em Let $I$ be a monomial ideal in $R$.  We say $I$ is a {\em linear quotient ideal}, if for some order 
$f_1,\ldots, f_m$ of the elements in $G(I)$ the colon ideal  $(f_1,\ldots ,f_{i-1}):f_i$ is generated by 
monomials of degree 1 for each $i\in[m]$, and all $f_i$ have the same degree. }
\end{Definition}

\begin{Definition}
\label{linear q-tree}
{\em Let $\Delta$ be a tree. If its facet ideal $I$  is a linear quotient ideal , then we call 
$\Delta$ a {\em linear quotient tree}. If $I$ has linear resolution, then we call $\Delta$ a 
{\em linear tree}.}
\end{Definition}

\begin{Proposition}
\label{linear resolution}
Let $\Delta $ be a tree, $I$ its facet ideal. 
\begin{enumerate}
\item[(i)] The following statements are equivalent:
  \begin{enumerate}
  \item[(a)] $\Delta$ is a linear quotient tree.
  \item[(b)]  $\Delta$ is a linear tree.
  \end{enumerate}
\item[(ii)] If $\Delta$ satisfies the equivalent conditions in {\em (i)}, then $\Delta$ is pure and connected in %%@
codimension
$1$.
\end{enumerate}
\end{Proposition}

\begin{proof}
(i) $(a)\Rightarrow(b)$: By definition, we have $\Delta$ is pure. Let $d-1$ be 
the dimension of $\Delta$. We prove the assertion  by induction on the number of facets $m$.

The case $m=1$ is trivial.  Suppose $m>1$. Let $F_1,\ldots,F_m$ be the facets of $\Delta$ and 
$J=(f_1,\ldots ,f_{m-1})$ such that $J:I$ can be generated by monomials of degree $1$. 
By induction hypothesis, $R/J$ has a linear resolution. Hence 
\begin{eqnarray}
\label{*}
\Tor_i^R(K,R/J)_j=0 \text{ for } j\neq i+d \text{ and all } i>0.      
\end{eqnarray}
Since $J:I$ is generated by monomials of degree 1, we have 
\begin{eqnarray}
\label{***}
\Tor_i^R(K,J:I)_j=0 \text{ for } j\neq i+1\text{ and all } i>0. 
\end{eqnarray}
From the exact sequence
\[
\begin{CD}
0@>>> R/(J:I)(-d)@>>> R/J@>>> R/I@>>> 0, 
\end{CD}
\]
we have the long exact sequence
\[
\begin{CD}
\cdots @>>> \Tor_i^R(K,R/(J:I)(-d))@>>> \Tor_i^R(K,R/J)@>>> \Tor_i^R(K,R/I)\\
@>>> \Tor_{i-1}^R(K,R/(J:I)(-d))@>>> \cdots,  
\end{CD}
\]
By using (\ref{*}) and (\ref{***}), this long exact sequence implies that  $\Tor_i^R(K,R/I)_j=0$ for 
$j\neq i+d$ and all $i>0$, 
so $I$ has linear resolution, i.e. $\Delta$ is a linear tree.

$(b)\Rightarrow (a)$: It is clear that if $\Delta$ is not pure, then $I$ has no linear resolution. We may assume 
$\Delta$ is a pure tree of dimension $d-1$. Suppose $I$ is not a linear quotient ideal. Let $F_1,\ldots,F_m$ be a leaf %%@
order. Then  $L=(f_1,\ldots ,f_{k-1}):f_k$ is not  generated by monomials of degree 1 for some $k\in\{1,\ldots,m\}$, %%@
and hence $b_{1,1+j}(R/L)\neq 0$ for some $j>1$. Let $I'=(f_1,\ldots,f_k)$ and $J'=(f_1,\ldots,f_{k-1})$. 
By Theorem \ref{basic} we have the exact sequence 
\[
\begin{CD}
0@>>> \Tor_2^R(K,R/J')@>>> \Tor_2^R(K,R/I')@>>> \Tor_1^R(K, R/L(-d))@>>> 0, 
\end{CD}
\]
which implies that  $b_{2,2+j+d}(R/I')\neq 0$, so $I'$ has no 
linear resolution since $I'$ is generated in degree $d$. By Corollary \ref{injective}, $I$ has no linear resolution, a %%@
contradiction.

(ii) It is clear that $\Delta$ must be pure. Let $F_1,\ldots,F_m$ be the facets of $\Delta$ such that
$(f_1,\ldots,f_{k-1}):f_k$ is generated by monomials of  degree 1 for $k=1,\ldots, m$. We prove that $\Delta$ is %%@
connected in
codimension $1$ by induction on $m$. The case $m=1$ is trivial. Assume $m>1$, since 
$(f_1,\ldots,f_{m-1})$ is a linear quotient ideal, by induction hypothesis, $\langle F_1,\ldots,F_{m-
1}\rangle$ is connected in codimension $1$. To show $\Delta$ is connected in codimension $1$, we only need 
to show that for any facet $F_i$, with $i<m$, there exists a proper chain between $F_i$ and $F_m$. Since 
$(f_1,\ldots,f_{m-1}):f_m$ is generated by monomials of degree $1$, we have that all the facets of 
$\langle F_m\rangle\cap\langle F_1,\ldots,F_{m-1}\rangle$ are of dimension $d-1$. Hence there exists an integer  
$j< m$ such that $\dim(F_j\cap F_m)=d-1$. Since $F_i$ and $F_j$ both are facets of the tree 
$\langle F_1,\ldots,F_{m-1}\rangle$, there exists a proper chain $F_i=F_{i_0},\ldots,F_{i_l}=F_j$ between 
$F_i$ and $F_j$. Hence $F_i=F_{i_0},\ldots,F_{i_l}=F_j,F_m$ is a proper chain between $F_i$ and $F_m$.
\end{proof}

By Proposition \ref {linear part} and Proposition \ref {linear resolution}  the 
Betti numbers of a linear tree  can now be described as follows:

\begin {Corollary}
\label {Betti numbers of lqt}
Let $\Delta$ be a $d$-dimensional linear  tree with $m$ facets, $I$ its facet ideal. Then 
\[
 b_i(R/I)= \left\{\begin{array}{ll}
               m,  &\text{if $i=1$,}\\
               \sum_{G}\binom{m(G)}{i}, &\text{if $i\geq 2$,}
                           \end{array}
                           \right.
\] 
where the sum is taken over all $(d-1)$-dimensional faces $G$  of $\Delta$,  and 
$m(G)=|\{F\in{\mathcal F}(\Delta)\: G\subset F\}|$. 
\end{Corollary}

Later in this section, we will classify  all linear trees of a given  dimension. For this, we need some 
preparation.

\begin{Lemma}
\label{part}
Let $\Delta$ be a linear tree, $\Gamma$  a subcomplex of $\Delta$ which is connected in codimension $1$. 
Then $\Gamma$ is a linear tree.  
\end{Lemma}

\begin{proof} 
It is clear that $\Gamma$ is again a pure tree. We may assume 
$\Gamma\neq\Delta$. We claim there exists an order of the facets $F_1,\ldots,F_l$ of 
$\langle {\mathcal F}(\Delta)\setminus{\mathcal F}(\Gamma)\rangle$ such that 
$\langle{\mathcal F}(\Gamma),F_1,\ldots,F_i\rangle$ 
is connected in codimension $1$, $i=1\ldots,l$. In fact, let $F\in{\mathcal F}(\Gamma)$ and 
$G\in{{\mathcal F}(\Delta)\setminus {\mathcal F}(\Gamma)}$ be any two facets. Since $\Delta$ is connected in
codimension $1$, there exists a unique irredundant proper chain from $F$ to $G$. Let $F_1$ be the first 
facet in this chain which does not belong to $\Gamma$. Then it is obvious that 
$\langle {\mathcal F}(\Gamma),F_1\rangle$ is connected in codimension $1$. 
The claim  follows by induction on $|{\mathcal F}(\Delta)\setminus {\mathcal F}(\Gamma)|$.

By Corollary \ref{leaf of original}, $F_i$ is a leaf of $\langle{\mathcal F}(\Gamma),F_1,\ldots,F_i\rangle$ 
for $i=1,\ldots,l$.  Let $I$ and $J$  be the facet ideals of 
$\Delta$ and $\Gamma$, respectively. By Corollary \ref{injective}, $b_{i,i+j}(J)\leq b_{i,i+j}(I)$ for any $i$
and $j$. Since $I$ has linear resolution, this implies that $J$ has linear resolution.
\end{proof}

\begin{Lemma}
\label{unique order}
Let $\Delta$ be a linear tree, $F$ and $G$ any two facets of $\Delta$. Let $F=F_0,\ldots,F_m=G$ be the 
irredundant proper chain between $F$ and $G$. Then 
$(f_0,\ldots,\\
 f_{l-1}):f_l$ is  generated by monomials of degree $1$, $l=0,\ldots,m$.
\end{Lemma}

\begin{proof}
Since $F_0,\ldots,F_m$ is an irredundant proper chain, $\langle F_0,\ldots,F_i\rangle$ is a linear tree 
for all $i$, see Lemma \ref{part}. Assume there exists an $l$ such that $(f_0,\ldots, f_{l-1}):f_l$ is 
not generated by monomials of degree $1$. Since $F_l$ is a leaf of $\langle F_0,\ldots,F_l\rangle$, it 
follows from Theorem \ref{basic} that $\langle F_0,\ldots,F_l\rangle$ is not a linear tree, 
a contradiction.
\end{proof}

\begin{Proposition}
\label{unique distance}
Let $\Delta$ be a pure  tree of dimension $d$, $F$ and $G$ any  two facets with $\dim(F\cap G)=d-k$, 
for some $k\in [d+1]$. Then 
\begin{enumerate}
\item[(i)] $\dist(F,G)\geq k$;
\item[(ii)] $\dist(F,G)=k$, if $\Delta$ is a linear tree.
\end{enumerate}
\end{Proposition}

\begin{proof}
(i) is obvious. Now let $\Delta$ be a linear tree, and suppose that $\dist(F,G)>k$. 
Let $F=F_0,\ldots ,F_l=G$ be the irredundant proper chain between $F$ and $G$, where $l>k$. 
Let $H=F\sect G$. By Proposition \ref{link}, $H\subset F_k$ for $k=0,\ldots ,l$. 

Let $\{x_{i}\}=F_{i}\setminus F_{i+1}$ for $i=0,\ldots, l-1$. We claim that 
$\{x_0,\ldots, x_{l-1}\}\subset F_0$, and that the  elements $x_i$ are pairwise distinct. 

Assume $x_j\notin F_0$ for some $j=0,\ldots, l-1$. Since $F_0,\ldots, F_j$ is an irredundant proper 
chain, it follows that $F_k\sect F_{j+1}$ is a proper subset of $F_{j}\sect F_{j+1}$ for $k<j$. This 
implies that $|F_k\setminus F_{j+1}|>1$  for all $k<j$, while $F_{j}\setminus F_{j+1}=\{x_j\}$. On the 
other hand, $(f_0,\ldots, f_j):f_{j+1}$ is generated by monomials of degree 1. This implies that 
$x_j\in F_k$ for all $k\leq j$. In particular, $x_j\in F_0$, a contradiction. Since $F_i,\ldots,F_l$ is an
irredundant proper chain, $F_i$ is a leaf of $\langle F_i,\ldots,F_l\rangle$ for all $i\in\{0,\ldots,l-1\}$. 
Hence $x_i\notin F_k$ for all $k>i$. So the $x_i$ are pairwise distinct,  and $x_i\notin H$ for 
$i=0,\ldots, l-1$.

So we have $H\cup \{x_0,\ldots, x_{l-1}\}\subseteq F_0$. Hence $|F_0|\geq d-k+1+l>d+1$, a contradiction.
\end{proof}

\begin{Definition}
\label{intersection property}
{\em Let $\Delta$ be a $d$-dimensional pure tree and connected in codimension $1$. If for any two facets $F$ 
and $G$ with $\dim(F\cap G)=d-k$, $k=1,\ldots ,d+1$, we have $\dist(F,G)=k$, then we say $\Delta$ has the 
{\em intersection property}.}
\end{Definition}

\begin{Remark}
\label{diameter of Delta}
{\em Let $\Delta$ be a $d$-dimensional tree with intersection property, and $l$ the diameter of $\Delta$. 
Then 
\begin{enumerate}
\item[(i)] $l\leq d+1$, and 
\item[(ii)] for any irredundant proper chain  ${\mathcal C}$  in $\Delta$, and  any  face $H$ in $\Gamma$ of 
dimension $d-k$, where $\Gamma$ is the simplicial complex generated by ${\mathcal C}$, one has that $H$  is 
contained in at most $k+1$ facets of $\Gamma$.
\end{enumerate}

In fact, it is clear that for any two facets $F$ and $G$ of $\Delta$, $\dist(F,G)\leq d+1$. 
Hence $l\leq d+1$. 

Assume $H$ is contained in more than $k+1$ facets of $\Gamma$. Since $\Gamma$ is generated by the irredundant 
proper chain ${\mathcal C}$, there exist two facets $F$ and $G$ of $\Gamma$ such that $H\subseteq F\cap G$ and 
$\dist(F,G)>k$. But $\dim(F\cap G)\geq \dim H=d-k$, contradicting Proposition \ref{unique distance}.}
\end{Remark}

\begin{Proposition}
\label{add facet}
Let $\Delta$ be a linear  tree of dimension $d$, and $G$ an adjacent face. Let 
$\Gamma=\langle {\mathcal F}(\Delta), F\rangle$, where $F$ is a facet of dimension $d$ and 
$\langle F\rangle\cap \Delta=\langle G\rangle$. Then $\Gamma$ is a linear tree. 
\end{Proposition}

\begin{proof}
By Proposition \ref{linear resolution}, we have  $\Delta=\langle F_1,\ldots ,F_m\rangle$ such that $(f_1,\ldots %%@
,f_{i-1}):f_i$ is generated by 
monomials of degree $1$, $i=1,\ldots,m$. Let $F_{i_1},\ldots ,F_{i_l}$ be all the facets of $\Delta$ 
which contains $G$, and $x_{i_j}=F_{i_j}\setminus F$ for $j=1,\ldots,l$, where $l>1$. We prove that
$(f_1,\ldots ,f_m):f=(x_{i_1},\ldots ,x_{i_l})$ (which implies that $\Gamma$ is also a linear tree).

It is clear that $(x_{i_1},\ldots ,x_{i_l})\subseteq (f_1,\ldots ,f_m):f$. In order to prove the converse inclusion , %%@
we first notice that there exists no facet $F_p$ of $\Delta$, such that $F_p\cap F_{i_j}=\emptyset$ for all %%@
$j=1,\ldots ,l$. 
Otherwise by Proposition \ref {unique distance}, $\dist(F_{i_j},F_p)=d+1$ for all $j=1,\ldots,l$. Since $l>1$, this %%@
contradicts  Lemma \ref{unique chain}.

It remains to show that for any facet $F_p$ of $\Delta$ we have $F_p\cap \{ x_{i_1},\ldots ,x_{i_l}\}\neq \emptyset$.
Suppose there exists a facet $F_p$ such that $F_p\cap \{ x_{i_1},\ldots ,x_{i_l}\}=\emptyset$, 
then $p\neq i_j$, and hence we have $F_p\cap G=F_p\cap F_{i_j}\neq \emptyset$ for  $j=1,\ldots ,l$. Let 
$\dim(F_p\cap F_{i_j})=d-k$. Then by Proposition \ref{unique distance}, $\dist(F_p,F_{i_j})=k$ for 
$j=1,\ldots,l$. Again, since $l>1$,  this contradicts  Lemma \ref{unique chain}.
\end {proof}

Now we can show

\begin{Theorem}
\label{all linear trees}
Let $\Delta$ be a tree. Then the following are equivalent:
\begin{enumerate}
\item[(i)] $\Delta$ is a linear tree.
\item[(ii)] $\Delta$ has intersection property.
\end{enumerate}
\end{Theorem}
\begin{proof}
 (i)\implies (ii) follows from Proposition \ref{unique distance}.

(ii)\implies (i): We prove the assertion by induction on the number of facets $m$ of $\Delta$. The case 
$m=1$ is trivial. Assume $m>1$. Let $F$ be a leaf of $\Delta$. By induction hypothesis, 
$\langle {\mathcal F}(\Delta)\setminus \{F\}\rangle$ is a linear tree because it still satisfies the intersection %%@
property. 
Let $H=\langle F\rangle\cap \langle {\mathcal F}(\Delta)\setminus \{F\}\rangle$; if $|{\mathcal U}_{\Delta}(F)|>1$,  %%@
then $H$ is an adjacent face of $\langle {\mathcal F}(\Delta)\setminus \{F\}\rangle$. Hence $\Delta$ is a linear tree %%@
by Proposition \ref{add facet}.  If $|{\mathcal U}_{\Delta}(F)|=1$, let $\{F'\}={\mathcal U}_{\Delta}(F)$ and 
$\{x\}=F'\setminus F$.  

 We claim $x$ is contained in any facet of $\langle {\mathcal F}(\Delta)\setminus \{F\}\rangle$. Hence, since  
 $\langle {\mathcal F}(\Delta)\setminus \{F\}\rangle$ is a linear tree, $\Delta$ is a linear tree, too.

In order to prove the claim, consider  $G\in {\mathcal F}(\Delta)$, $G\neq F$, and let $E=F\sect G$ and assume that %%@
$\dim E=d-k$. Then 
$G=\{x_1,\ldots, x_k\}\cup E$ and $F=\{y_1,\ldots, y_k\}\cup E$, where all the elements in  %%@
$\{x_1,\ldots,x_k,y_1,\ldots,y_k\}$ are pairwise  distinct. Since 
$\Delta$ has intersection property, $\Delta$ is pure and connected in codimension $1$, and $\dist(F,G)=k$. 
Hence there exists an irredundant proper chain $G=F_0,F_1,\ldots F_k=F$ between $G$ and $F$. Since $F$ is 
a leaf of $\Delta$ and $\{F'\}={\mathcal U}_{\Delta}(F)$, we have $F_{k-1}=F'$. Since 
$|F_i\setminus F_{i+1}|=1$ for all $i$, we may assume $F_i=\{y_1,\ldots,y_i,x_{i+1},\ldots, x_k\}\cup E$ for 
$i=1,\ldots,k$. Hence $F_{k-1}=\{y_1,\ldots,y_{k-1}, x_k\}\cup E$. But on the other hand, 
$F_{k-1}=F'=\{y_1,\ldots,y_{k-1},x\}\cup E$. Hence $x=x_k\in G$. 
\end{proof}

\begin{Remark}
\label{on stanley reisner rings}
{\em In the case that $\Delta$ is a 1-dimensional tree, the intersection property is equivalent to the condition that 
the distance between any two edges in $\Delta$ is at most 2, and this is equivalent to say that the complement  
$\bar{\Delta}$ of the graph $\Delta$ is a triangulated graph. This coincides with the result of Fr\"oberg in \cite %%@
{F}.}
\end{Remark}

\section{ The alternating sum property  of facet ideals}

In this section we show that for a special  class of  facet ideals $I$ the Betti numbers have the property   that %%@
$\sum_i(-1)^ib_{i,i+j}(R/I)=0$ for all $j>d$, where $d$ is the least  degree of the generators.  These class of ideals %%@
include facet ideals  of trees (not necessary pure) which are connected in codimension 1.

\begin{Definition}
\label{sum}
{\em  Let $I$ be a monomial ideal in $R$ with  $G(I)=\{f_1\ldots,f_m\}$ and  $d=\min\{\deg(f_i)\: i=1,\ldots,m\}$. We %%@
say that $I$ has the {\em alternating sum property}, if
\[
\sum_{i\geq 1}{(-1)^i b_{i,i+j}(R/I)}=\left\{ \begin{array}{ll}
                                              -1,  &\text{for $j=d$,}\\
                                               0,  &\text{for $j>d$.}
                                           \end{array}
                                              \right.
\]
}
\end{Definition}

To proof the main theorem of this section, we need the following fact:

\begin{Lemma}
\label{single variable}
Let $I$ be a monomial ideal in $R$. Suppose $G(I)$ contains a monomial of degree $1$. Then 
$\sum_i{(-1)^ib_{i,i+j}(R/I)}=0$ for all $j$.
\end{Lemma}

\begin{proof}
Let $G(I)=\{m_1,\ldots,m_l,x\}$, $J=(m_1,\ldots,m_l)$. Then $x$ does not divide $m_j$ for $j=1,\ldots,l$, and 
$J:I=J$. By Theorem \ref{basic} we have for $i>0$  
\[
\begin{CD}
b_{i,i+j}(R/I)&=&b_{i,i+j}(R/J)+b_{i-1,i+j}(R/(J:I)(-1))\\
&=&b_{i,i+j}(R/J)+b_{i-1,i-1+j}(R/(J:I))\\
&=&b_{i,i+j}(R/J)+b_{i-1,i-1+j}(R/J).
\end{CD}
\]
From this it follows that $\sum_i{(-1)^ib_{i,i+j}(R/I)}=0$. 
\end{proof}

\begin{Remark}
{\em With the same arguments as in the proof of Lemma \ref{single variable} one can show more generally: Let $J$ be a %%@
graded ideal in $R$, $I=(J,f)$, where $\deg(f)=1$. If $f$ is regular on $R/J$,
 then  $\sum_i{(-1)^ib_{i,i+j}(R/I)}=0$ for all $j$.}
\end{Remark}

\begin{Proposition}
\label{main}
Let $\Delta$ be a simplicial complex with facet ideal $I$. 
If there exists an order of the facets  $F_1,\ldots,F_m$  of $\Delta$ 
such that for each $i=2,\ldots,m$, $F_i\setminus{\bigcup_{j<i} F_j}\neq\emptyset$, and there exists $j<i$ such that 
$|F_j\setminus F_i|=1$. Then $I$ has the alternating sum property.
\end{Proposition}

\begin{proof}
We prove this proposition  by induction on $m$. The case $m=1$ is trivial. Let $d=\min\{\deg(f_i)\:i=1,\ldots,m-1\}$, 
$d'=\deg(f_m)$, and $J=(f_1,\ldots,f_{m-1})$. 
Since $|F_j\setminus F_m|=1$ for some $j<m$ it follows that  $d'\geq d$, and that  $G(J:I)$ contains at least one %%@
monomial of degree $1$. By 
Lemma \ref{single variable}, 
\begin{eqnarray}
\label{s}
\sum_i{(-1)^ib_{i,i+j}(R/(J:I))}=0 \text{  for any  } j.
\end{eqnarray}
On the other hand by Theorem \ref{basic}, we have 
\begin{eqnarray}
\label{b}
b_{i,i+j}(R/I)=b_{i,i+j}(R/J)+b_{i-1,i+j-d'}(R/(J:I)), 
\end{eqnarray}
for $i>0$, since  $F_m\setminus{\bigcup_{j<m} F_j}\neq\emptyset$. By induction hypothesis $J$ has the alternating sum %%@
property.
Hence one sees that $I$ has the alternating sum property by using (\ref{s}) and (\ref{b}). 
\end{proof}

\begin{Corollary}
\label{pure quasi}
Let $\Delta$ be a pure quasi-tree connected in codimension 1 with facet ideal $I$.  Then $I$ has the alternating sum %%@
property.
\end{Corollary}

\begin{proof}
 Since $\Delta$ is a quasi-tree, there exists a leaf order of facets $F_1,\ldots,F_m$. The assertion 
follows from Proposition \ref{main} immediately.
\end{proof}

The next result shows that in Corollary \ref{pure quasi} we can skip the assumption that $\Delta$ is 
pure if we assume that $\Delta$ is a tree.

\begin{Theorem}
\label{not pure tree}
Let $\Delta$ be a tree connected in codimension $1$ with facet ideal $I$. Then $I$ has the alternating 
sum property.
\end{Theorem}

\begin{proof}
 We prove the assertion by induction on the number of facets $m$. The case $m=1$ is trivial. Assume 
$m>1$.  Let $d=\dim\Delta$. There are two cases. 

Case 1. There exists only one facet $F$ of dimension $d$. Then $F$ must be a leaf. Otherwise, there exist 
two facets $G_1$, $G_2$ such that $F\cap G_1\not\subseteq F\cap G_2$ and 
$F\cap G_2\not\subseteq F\cap G_1$. Since $\Delta$ is connected in codimension 1 and $\dim G_i<d$, 
$i=1,2$, there exists a chain ${\mathcal C}$ between $G_1$ and $G_2$ which does not include $F$. Then 
the simplicial subcomplex $\Gamma$ whose facets are the elements of ${\mathcal C}$ and $F$ has no leaf, 
a contradiction. 

We choose a  $G\in{\mathcal U}_{\Delta}(F)$ (see Definition \ref{leaf}) of maximal  dimension. Since 
$\Delta$ is connected in codimension $1$, we have $\dim G=\dim\langle{\mathcal F}(\Delta)\setminus \{F\}\rangle$ 
and $\dim (F\cap G)=\dim G-1$, i.e. $|G\setminus F|=1$. 
Since $F$ is a leaf, $\langle{\mathcal F}(\Delta)\setminus \{F\}\rangle$ is a tree with $m-1$ facets 
which is connected in codimension $1$. By induction hypothesis there exists a leaf order of facets 
$F_1,\ldots,F_{m-1}$ such that for each $i=2,\ldots,m-1$,    
$F_i\setminus{\bigcup_{j<i} F_j}\neq\emptyset$, and there exists $j<i$ such that 
$|F_j\setminus F_i|=1$. Let $F=F_m$. We see that $F_1,\ldots,F_m$ satisfy the conditions of 
Theorem \ref{main} in this order.

Case 2. There exist more than one facets of dimension $d$. Let $G_1,\ldots,G_s$ be all of these facets, 
where $s>1$. Then for any $i$ and $j$, the facets in any proper chain between $G_i$ and $G_j$  are all 
of dimension $d$, and hence belong to $\{ G_1,\ldots,G_s\}$. Therefore 
$\Sigma=\langle G_1,\ldots,G_s\rangle$ is pure tree and  connected in codimension 1.  
By Proposition \ref{two leaves}, $\langle G_1,\ldots,G_s \rangle$ has at least two leaves. 

We claim that at least one of the leaves of  $\Sigma$  is a leaf of $\Delta$. Suppose this is not the 
case. We take any two leaves of $\Sigma$, say $G_i$ and $G_j$ with free vertex $x_i$ and $x_j$, 
respectively. Since $G_i$ and $G_j$ are not leaves in $\Delta$ there exist elements  
$F, F'\in {\mathcal F}(\Delta)\setminus {\mathcal F}(\Sigma)$ with $x_i\in F$ and $x_j\in F'$. Let 
${\mathcal C}$ be a chain between $F$ and $F'$.  Since $\dim F<d$ and $\dim F'<d$, all elements of this 
chain do not belong to ${\mathcal F}(\Sigma)$. On the other hand, let ${\mathcal C}'$ be a proper chain 
between $G_i$ and $G_j$, then all elements of the chain belong to ${\mathcal F}(\Sigma)$, because 
$\dim G_i=\dim G_j=d$. Then the simplicial complex generated by the elements of these two chains has no 
leaf, a contradiction. 

We may assume that  $G_i$ is a leaf of $\Delta$.  Removing $G_i$ from $\Delta$ yields a tree which is 
again connected in codimension 1, and we may proceed as in case 1. 
\end{proof}

\begin{Corollary}
\label{gragh tree}
Let $\Delta$ be a $1$-dimensional tree with facet ideal $I$. Then $I$ has the alternating sum property.
\end{Corollary}

\begin{proof}
It is clear that $\Delta$ is connected in codimension $1$. The result follows from Theorem \ref{not pure
tree}.
\end{proof}

\end{document}